\documentclass[reqno]{amsart}
\usepackage{amsmath}
\usepackage{amsfonts}
\usepackage{amssymb}
\usepackage{amsthm}
\newtheorem{theorem}{Theorem}[section]
\newtheorem{lemma}[theorem]{Lemma}
\newtheorem{proposition}[theorem]{Proposition}
\newtheorem{corollary}[theorem]{Corollary}
\newtheorem{notation}[theorem]{Notation}
\newtheorem{definition}[theorem]{Definition}
\newtheorem{example}[theorem]{Example}
\newtheorem{remark}[theorem]{Remark}

\numberwithin{equation}{section}

\begin{document}
\title[Arithmetic structures on RM noncommutative tori]{Arithmetic structures on noncommutative tori with real multiplication}
\author{Jorge Plazas}
\address{ Max Planck Institute for Mathematics,
Vivatsgasse 7,
Bonn 53111,
Germany}
\email{plazas@mpim-bonn.mpg.de}

\date{}

\maketitle

\begin{abstract}
We study the homogeneous coordinate rings of real multiplication noncommutative tori as defined in \cite{Po}. Our aim is to understand how these rings give rise to an arithmetic structure on the
noncommutative torus. We start by giving an explicit presentation of these rings in terms of their natural generators. This presentation gives us information about the rationality of the rings and is the starting point for the analysis of the corresponding geometric data. Linear bases for these rings are also discussed. The basic role played by theta functions in these computations gives relations with elliptic curves at the level of modular curves. We use the modularity of these coordinate rings to obtain algebras which do no depend on the choice of a  complex structure on the noncommutative torus. These algebras are obtained by an averaging process over (limiting) modular symbols. 
\end{abstract}
\setcounter{section}{-1}

\section{Introduction}
\label{intro}

Noncommutative tori are standard prototypes of noncommutative spaces. Since the early stages of noncommutative geometry these spaces have been central examples
arising naturally in various context. In \cite{manin} Manin proposed the use of noncommutative tori as a geometric framework for the study of abelian class field theory of real  quadratic fields. This is the so called ``real multiplication program". The main idea is that noncommutative tori may play a role in the study of real quadratic fields analogous to the role played by elliptic curves in the study of imaginary quadratic fields. The relation between the 
endomorphism rings of noncommutative tori and orders in real quadratic fields is a good evidence supporting this point of view. If a noncommutative torus admits nontrivial Morita autoequivalences then the ring of such autoequivalences is an order in a real quadratic field. A noncommutative torus having this property is called a real multiplication noncommutative torus. 

The fact that noncommutative geometry may be relevant for addressing 
questions in number theory is also supported by various results and relations that have emerged in the last years. In particular explicit class field theory of $\mathbb{Q}$ (Kronecker-Weber theorem) and explicit class field theory of imaginary quadratic fields (complex multiplication) can both be recovered from the dynamics of certain quantum statistical mechanical systems (\cite{BostConnes,CMR1,CMR2}). The existence of quantum statistical mechanical systems with rich arithmetical properties opens a new approach to the study of explicit class field theory using the tools of quantum statistical mechanics. The first case for which there is not yet a complete solution to the explicit class field theory problem is the case of real quadratic fields, $K = \mathbb{Q}(\sqrt{D})$, where $D \in \mathbb{N}^{+}$ is a square free positive integer. 

Noncommutative tori are, a priori, analytical objects. In order to achieve arithmetical applications it is important to find appropriate algebraic structures underlying these spaces. Rings admitting models algebraic over the corresponding base fields have proved to be essential for the analysis of quantum statistical mechanical systems of arithmetic nature. 

In \cite{Po} Polishchuk defined homogeneous coordinate rings for real multiplication noncommutative tori endowed with a complex structure. These rings seem to be good candidates for the applications described above. Our aim is to understand in which sense these rings provide an arithmetic structure on noncommutative tori. Starting from the explicit formulae defining these rings we study their rationality properties and their dependence on the complex structure. Some examples related to particular quadratic fields are treated. 

Homogeneous coordinate rings for real multiplication noncommutative tori can be viewed  as a family of algebras varying with the parameter defining the complex structure on the noncommutative torus. We give a presentation of these rings in terms of modular forms and use this modularity to define rings which do not depend on the choice of a complex structure. This is done by an averaging process over (limiting) modular symbols.
The fact that modular forms play an essential role in our analysis points to deep relations with the quantum thermodynamical system introduced by Connes and Marcolli in \cite{ConnesMarcolli}. This quantum mechanical system recovers the class field theory of the modular field and provides a two dimensional analog of the dynamical system corresponding to the Kronecker-Weber theorem.

In the case of real quadratic fields explicit class field theory is conjecturally given in terms of special values of $L$-functions, this is the content of Stark's famous conjectures \cite{Stark}. In order to apply our results on noncommutative tori in this direction using the tools of quantum statistical mechanics we still need to find $C^{*}$-completions for the homogeneous coordinate rings of real multiplication noncommutative tori. As a preliminary step in this direction we describe the geometric data associated to the homogeneous coordinate rings of noncommutative tori. We expect that at a future stage the use of the techniques developed in \cite{CoDV} will make it possible to obtain suitable $C^{*}$-completions.
In a different but related perspective arithmetic structures on noncommutative tori have been recently studied by Vlasenko in \cite{Masha} where the theory of rings of quantum theta functions is developed.

I want to thank Matilde Marcolli for her continued guidance through this project. I am also grateful many people at MPI and UniBonn for very helpful discussions. Special thanks go to Nikolai Durov, Snigdhayan Mahanta and Eugene Ha. I thank the Max-Planck-Institut f\"{u}r Mathematik  for its support and hospitality. 


\section{Basic setting}
\label{sec1}
In this section we recall the basic facts about noncommutative tori that will be needed in the sequel. Also, we recall the definition of homogeneous coordinate rings of real multiplication noncommutative tori given in \cite{Po}.

Noncommutative  tori are defined by their algebras of 
smooth functions which are noncommutative deformations of  $C^\infty (\mathbb{T})$, 
the algebra of smooth functions on the two dimensional torus $\mathbb{T}=S^1\times S^1$. 
Given  $\theta \in \mathbb{R}$ we define $\mathcal{A}_\theta$, the {\it algebra of smooth functions on the noncommutative torus $\mathbb{T}_\theta$}, as the algebra of formal power series in two unitaries $U$ and $V$ with rapidly decreasing coefficients and multiplication given by the relation $UV=e^{2\pi i \theta}VU$:
\begin{eqnarray*}
\mathcal{A}_\theta &=& C^\infty (\mathbb{T}_\theta) \\
 &=& \{ a =
\sum_{n,m \in \mathbb{Z}} a_{n,m}U^nV^m  \, | \, \{a_{n,m} \}
\in \mathcal{S}(\mathbb{Z}^2) \} \nonumber
\end{eqnarray*}

Let $\theta\in\mathbb{R}$ be irrational. Following \cite{Connes1} we define, for any pair $c,d\in \mathbb{Z}$, $c>0$, a right $\mathcal{A}_\theta$-module $E_{d,c}(\theta)$ given by the following action of $\mathcal{A}_\theta $ on the Schwartz space $\mathcal{S} (\mathbb{R} \times \mathbb{Z}/ c\mathbb{Z}) = \mathcal{S}(\mathbb{R})^c$: 
\begin{eqnarray}
(fU)(x,\alpha)&=& f(x- \frac{c\theta+d}{c},\alpha -1) \\
(fV)(x,\alpha)&=& \exp (2 \pi i (x- \frac{\alpha d}{c})) f(x,\alpha)
\end{eqnarray}
The right $\mathcal{A}_\theta$-module $E_{d,c}(\theta)$ is projective and of finite type.
If $c$ and $d$ are relatively prime we say that $E_{d,c}(\theta)$ is a 
{\it basic $\mathcal{A}_\theta $-module}. Being this the case the pair $d,c$ can be completed to a matrix  
\begin{eqnarray}
\label{lag}
g= \left(
\begin{array}{cc}
  a & b \\
  c & d \\
\end{array}
\right)\in SL_2 (\mathbb{Z})
\end{eqnarray}
In this case we write $E_g(\theta)$  for the module $E_{d,c}(\theta)$. By
definition the degree of $E_g(\theta)$ is taken to be $c$. We also define the
degree of a matrix $g \in SL_2 (\mathbb{Z})$, given as above, by $\deg(g)=c$. The study of the differential geometry of $\mathbb{T}_\theta$ based on canonical
derivations on $\mathcal{A}_{\theta}$ and the corresponding compatible connections on the above modules 
carried out in \cite{Connes1} and further developed in \cite{CoRiff} is the
starting point of the more recent algebro-geometric setup we describe below.  

Let $SL_2 (\mathbb{Z})$ act on $\mathbb{R}$ by fractional linear transformations. Let $g\in SL_2 (\mathbb{Z})$ be as above and denote by $U'$ and $V'$ two generating unitaries of the algebra $ \mathcal{A}_{g \theta}$. We can define a left action of the algebra $ \mathcal{A}_{g \theta}$ on $E_g$ by: 
\begin{eqnarray}
(U'f)(x,\alpha)&=& f \left(x- \frac{1}{c},\alpha -a \right) \\
(V'f)(x,\alpha)&=& \exp (2 \pi i (\frac{x}{c\theta +d}-\frac{\alpha}{c})) f(x,\alpha)
\end{eqnarray}
This action gives an identification: 
\begin{eqnarray}
End_{\mathcal{A}_\theta}(E_g(\theta)) \simeq \mathcal{A}_{g\theta}.
\end{eqnarray}

Note in particular that if $g\theta =\theta$ then $E_g(\theta)$ has the 
structure of a $\mathcal{A}_\theta$-bimodule. Therefore in this case the algebra $\mathcal{A}_\theta$ has nontrivial Morita auto equivalences. An irrational number $\theta\in\mathbb{R}\setminus \mathbb{Q}$ is a fixed point of a fractional  linear transformation $g \in SL_2 (\mathbb{Z})$ if and only if it generates a quadratic extension of $\mathbb{Q}$. In this case we say that the noncommutative torus $\mathbb{T}_\theta$ with algebra of smooth functions $\mathcal{A}_\theta$ is a {\it real multiplication noncommutative torus}. 

The tensor product of basic modules is again a basic module 
(c.f. \cite{CoRiff, PoS,DiS}). Given $g_{1}, g_{2} \in SL_2 (\mathbb{Z})$, there is a well defined pairing of right $\mathcal{A}_\theta$-modules:
\begin{eqnarray}
\label{producto1}
t_{g_{1}, g_{2}}: E_{g_{1}}(g_{2}\theta)\otimes_{\mathbb{C}}  E_{g_{2}} (\theta) \rightarrow 
E_{g_{1} g_{2}}(\theta)
\end{eqnarray}
This map gives rise to an isomorphism of $\mathcal{A}_{g_{1}g_{2}\theta} - \mathcal{A}_{\theta}$ bimodules: 
\begin{eqnarray}
E_{g_{1}}(g_{2}\theta)\otimes_{\mathcal{A}_{g_{2}\theta}}  E_{g_{2}} (\theta) \rightarrow 
E_{g_{1} g_{2}}(\theta).
\end{eqnarray}
In particular, if $g\theta = \theta$ one has an isomorphism 
\begin{eqnarray}
\underbrace{E_{g}(\theta)\otimes_{\mathcal{A}_ {\theta}} \dots
\otimes_{\mathcal{A}_ {\theta}}  E_{g} (\theta)}_{n} 
\simeq  E_{g^{n}} (\theta) .
\end{eqnarray}

Given $\tau \in \mathbb{C}$ with $Im(\tau) < 0$ we endow the noncommutative torus with a complex structure given by the derivation 
$\delta_\tau:\mathcal{A}_\theta \rightarrow \mathcal{A}_\theta$:
\begin{eqnarray}
\delta_\tau: \sum_{n,m \in \mathbb{Z}} a_{n,m}U^nV^m \mapsto (2 \pi \imath )\sum_{n,m \in \mathbb{Z}} (n\tau + m)a_{n,m}U^nV^m 
\end{eqnarray}
This derivation should be viewed as an analog of the operator $\bar{\partial}$ on a complex elliptic curve. We will denote by $\mathbb{T}_{\theta, \tau}$ the noncommutative torus $\mathbb{T}_\theta$ equipped with this complex structure.

A holomorphic structure on a  right $\mathcal{A}_\theta$-module $E$ is given by an operator $\bar{\nabla}:E\rightarrow E$ which is compatible with the complex structure $\delta_\tau$ in the sense that it satisfies the following Leibniz rule:
\begin{eqnarray}
\bar{\nabla}(ea)= \bar{\nabla}(e)a + e \delta_\tau (a), \quad e \in
E, a\in\mathcal{A}_\theta 
\end{eqnarray}
Given a holomorphic structure $\bar{\nabla}$ on a right $\mathcal{A}_\theta$-module $E$ the corresponding set of holomorphic sections is the space $H^0(\mathbb{T}_{\theta,\tau}, E_{\bar{\nabla}}):= Ker (\bar{\nabla})$.

On every basic module $E_{d,c}$ one can define a family of holomorphic structures $\{ \bar{\nabla}_z \}$ depending on a complex parameter $z\in \mathbb{C}$:
\begin{eqnarray}
\label{nabla}
\bar{\nabla}_z(f)=\frac{\partial f}{\partial x} + 2 \pi i \left( \frac{ d \tau}{c \theta + d}
x+z \right )f.
\end{eqnarray}
By definition a {\it standard holomorphic vector bundle on $\mathbb{T}_{\theta, \tau}$} is given by a basic module $E_{d,c}=E_{g}$ together with one of the holomorphic structures $\bar{\nabla}_z$.

The spaces of holomorphic sections of a standard holomorphic vector bundles on $\mathbb{T}_{\theta, \tau}$ are finite dimensional (c.f. \cite{PoS}, Section 2). If $c\theta + d > 0$ then $\dim H^0(E_{g}, \bar{\nabla}_0 ) = c$. On what follows we will consider the spaces of holomorphic sections corresponding to  $ \bar{\nabla}_0$:
\begin{eqnarray}
\label{H0}
\mathcal{H}_g := H^0(\mathbb{T}_{\theta,
\tau}, E_{g, \bar{\nabla}_0}).
\end{eqnarray}

A basis of $\mathcal{H}_g$ is given by the Schwartz functions:
\begin{eqnarray}
\label{bases1}
\varphi_{\alpha}(x,\beta) = \exp (-\frac{c \tau }{c\theta +d} \frac{x^2}{2})\delta^{\beta}_{\alpha}  \qquad \alpha =1,...,c.
\end{eqnarray}

The tensor product of holomorphic sections is again holomorphic. Using the above basis the product can be written in terms of the corresponding structure constants.

\begin{theorem}\emph{(\cite{PoS} Section 2)}
\label{poli1}
Suppose $g_{1}$ and $g_{2}$ have positive degree. Then $g_1g_2$ has positive degree and $t_{g_{1}, g_{2}}$ induces a well defined linear map
\begin{eqnarray}
\label{producto2}
t_{g_{1}, g_{2}}: \mathcal{H}_{g_{1}}(g_{2}\theta)\otimes_{\mathbb{C}}\mathcal{H}_{g_{2}} (\theta) \rightarrow \mathcal{H}_{g_{1} g_{2}}(\theta).
\end{eqnarray}
Let $g_1,g_2$ and $g_1g_2$ be given by  
$$
g_1= \left(
\begin{array}{cc}
  a_1 & b_1 \\
  c_1 & d_1 \\
\end{array}
\right), \quad 
g_2= \left(
\begin{array}{cc}
  a_2 & b_2 \\
  c_2 & d_2 \\
\end{array}
\right), \quad 
g_1 g_2 = \left(
\begin{array}{cc}
  a_{12} & b_{12} \\
  c_{12} & d_{12} \\
\end{array}
\right)
$$
and let $\{ \varphi_{\alpha} \}$, $\{\varphi^{'}_{\beta} \}$ and $ \{ \psi_{\gamma} \}$ be respectively the basis of $\mathcal{H}_{g_{1}}(g_{2}\theta)$, $\mathcal{H}_{g_{2}} (\theta)$ and $\mathcal{H}_{g_{1} g_{2}}(\theta)$ as given in \ref{bases1}. Then 
\begin{eqnarray}
t_{g_{1}, g_{2}}:  \varphi_{\alpha}\otimes \varphi^{'}_{\beta} \mapsto C^{\gamma}_{\alpha , \beta} \psi_{\gamma}
\end{eqnarray}
where for $\alpha= 1,...,c_1$, $\beta = 1,...,c_2$ and $\gamma = 1,..., c_{12}$ we have:
\begin{eqnarray}
\label{C1}
C^{\gamma}_{\alpha , \beta} = \sum_{m\in I_{g_1, g_2}(\alpha , \beta, \gamma)} 
\exp [\pi \imath \frac{- \tau m^{2}}{2  c_1  c_2 c_{12}} ]
\end{eqnarray}
with
\begin{eqnarray*}
I_{g_1, g_2}(\alpha , \beta, \gamma) = \{ n\in \mathbb{Z} &|&  n\equiv -c_1\gamma + c_{12}\alpha \mod  c_{12}c_1, \; \\
&&n\equiv c_{2}d_{12} \gamma - c_{12}d_2 \beta \mod c_{12}c_2 \}.
\end{eqnarray*}
\end{theorem}
\begin{notation}
Throughout we use the convention of summing over  repeated indexes.
\end{notation}


Assume now that $\theta \in \mathbb{R}$ is a quadratic irrationality. So there exist some $g \in SL_2(\mathbb{Z})$ with $g\theta= \theta$ and $\mathbb{T}_{\theta}$ has real multiplication. Fix a complex structure $\tau$ on $\mathbb{T}_{\theta}$. In the case 
$E=E_g(\theta)$ we can extend a holomorphic structure on $E_g$ to a holomorphic structure on the tensor powers $E_{g}^{\otimes n}$. Following \cite{Po} we define a {\it homogeneous coordinate ring for $\mathbb{T}_{\theta,\tau}$} by 
\begin{eqnarray}
\label{LAdefinicion}
B_g(\theta,\tau) &=& 
\bigoplus_{n\geq 0} H^0(\mathbb{T}_{\theta,\tau}, E_{\bar{\nabla}_0}^{\otimes n})\\
 &=& \bigoplus_{n\geq 0} \mathcal{H}_{g^{n}} \nonumber
\end{eqnarray}
One of the main results in \cite{Po} characterizes some structural properties of these algebras in terms of the matrix elements of $g$:
\begin{theorem}\emph{(\cite{Po} Theorem 3.5)}
\label{poli2}
Assume $g\in SL_2 (\mathbb{Z})$ has positive real eigenvalues 
\begin{itemize}
\item If $c \geq a+ d$ then $B_g(\theta,\tau)$ is generated over $\mathbb{C}$ by $\mathcal{H}_{g}$.
\item If $c \geq a+ d+1$ then $B_g(\theta,\tau)$ is a quadratic algebra.
\item If $c \geq a+ d+2$ then $B_g(\theta,\tau)$ is a Koszul algebra.
\end{itemize}
\end{theorem}

We will fix some notations and conventions for the rest of the paper. As above $g$ will denote a matrix
\begin{eqnarray*} 
g= \left(
\begin{array}{cc}
  a & b \\
  c & d \\
\end{array}
\right)\in SL_2 (\mathbb{Z})
\end{eqnarray*}
We will always assume that the following inequalities hold 
\begin{eqnarray}
\label{condiciones1}
Tr(g) &= a+d > 2  &(g \text{ is hyperbolic}) \\
\label{condiciones2}
c & \geq  Tr(g) +2 & =a+d + 2 
\end{eqnarray}
The first inequality implies that $g$ has positive real eigenvalues. By Theorem~\ref{poli2} the second inequality implies that $B_{g}(\theta,\tau)$ is a quadratic Koszul algebra generated in degree 1.

We denote by $\lambda^{+}$ and $\lambda^{-}$ the eigenvalues of $g$ with $ 0 < \lambda^{+}< 1$ and $1<\lambda^{-}$. It is important to note that $\lambda^{-}$ is a fundamental unit  for the quadratic extension it generates. We also take 
\begin{eqnarray}
\label{tetha1}
\theta = \frac{ \lambda^{+} - d}{c} \quad, \quad \theta ' = \frac{ \lambda^{-} - d}{c}
\end{eqnarray}
These are the fixed points of $g$. 
\begin{proposition}
Let $\alpha \in \mathbb{R}$ be a quadratic irrationality. Then there exist $g$ and $\theta$ satisfying \ref{condiciones1} and \ref{condiciones2} above and such that $\mathbb{Q}(\alpha)=\mathbb{Q}(\theta)$.  
\end{proposition}

\begin{proof}
Suppose $\alpha \in \mathbb{R}$ is a quadratic irrationality. Then there exist a hyperbolic element 
$$
h = \left(
\begin{array}{cc}
  a' & b' \\
  c' & d' \\
\end{array}
\right)\in SL_2 (\mathbb{Z})
$$
having $\alpha$ as one of its two fixed points. Being a fixed point of $h$, $\alpha$ satisfies the quadratic equation $c'\alpha ^2 + (d'-a') \alpha - b'=0$. Since $a'd'-b'c' =1$  we can write the discriminant of this equation as $D = (a'+d')^2 - 4 = Tr(h)^2 -4$. The quadratic irrationalities $\alpha$ and $\sqrt{D}$ generate the same field extension of $\mathbb{Q}$. $| Tr(h) | > 2 $ and we may assume $Tr(h) > 2 $ since multiplying $h$ by $-1$ does not change $D$. Define 
$$
g = \left(
\begin{array}{cc}
  Tr(h)+1 & -1 \\
  Tr(h)+2 & -1 \\
\end{array}
\right)
$$
Then $g\in SL_2 (\mathbb{Z})$ and $Tr(g)=Tr(h)$ so the fixed points of $g$ generate the same extension of $\mathbb{Q}$ than $\alpha$. By construction $g$ satisfies 
\ref{condiciones1} and \ref{condiciones2}. 
\end{proof}

\section{$B_{g}(\theta,\tau)$ in terms of generators and relations}
\label{sec2}
We want to describe $B_{g}(\theta,\tau)$ in terms of generators and relations. Let $\{ \varphi_{\alpha}| \alpha = 1,...,c \}$ be the basis for $\mathcal{H}_g$ and  
$\{ \psi_{\gamma}| \gamma = 1, ..., c(a+d) \}$ the corresponding basis of $\mathcal{H}_{g^2}$ as given in \ref{bases1}.

The multiplication map $m:\mathcal{H}_{g}\otimes \mathcal{H}_{g}\rightarrow \mathcal{H}_{g^2}$ of the algebra $B_{g}(\theta,\tau)$ is $t_{g,g}$. It is given in the above basis as in Theorem~\ref{poli1} by the structure constants:   
\begin{eqnarray}
\label{producto}
m: \varphi_{\alpha}\otimes \varphi_{\beta} \mapsto C^{\gamma}_{\alpha , \beta} \psi_{\gamma}
\end{eqnarray}

Our computation of the relations defining $B_{g}(\theta,\tau)$ is based on the following observation:
\begin{lemma}
\label{lemmafacil}
Denote by $\mathcal{T}$ the tensor algebra of $\mathcal{H}_g$. Then $B_{g}(\theta,\tau)$
is isomorphic to the quotient of $\mathcal{T}$ by the homogeneous ideal 
$\mathcal{R}$ generated by $Ker(m) \subset \mathcal{H}_g \otimes \mathcal{H}_g$. 
\end{lemma} 

Since $m: v^{\alpha , \beta}\varphi_{\alpha}\otimes \varphi_{\beta} \mapsto v^{\alpha , \beta} C^{\gamma}_{\alpha , \beta} \psi_{\gamma}$ we have that  $v^{\alpha , \beta}\varphi_{\alpha}\otimes \varphi_{\beta} \in \mathcal{H}_g \otimes \mathcal{H}_g$ belongs to $Ker(m)$ if and only if $v^{\alpha , \beta} C^{\gamma}_{\alpha , \beta} = 0 $ for all $\gamma = 1,..., c(a+d)$. Using the bases $\{ \varphi_{\alpha} \}$ and $\{ \psi_{\gamma}\}$  we identify $\mathcal{H}_{g}\otimes \mathcal{H}_{g}$ and $\mathcal{H}_{g^2}$ with $\mathbb{C}^{c^2}$ and $\mathbb{C}^{c(a+d)}$ respectively. Finding a set of defining relations of $B_{g}(\theta, \tau)$ for the generators $\{ \varphi_{\alpha}| \alpha = 1,...,c \}$ amounts to finding a basis for the kernel of the linear map $M:\mathbb{C}^{c^2} \rightarrow \mathbb{C}^{c (a+d)}$ 
given by the $C^{\gamma}_{\alpha , \beta}$.

\begin{lemma}
\label{lemmaconstantes}
The structure constant $ C^{\gamma}_{\alpha , \beta}$ is different from zero if and only 
if $ \alpha  \equiv d(\gamma - \beta) \mod c $.
\end{lemma}
\begin{proof}
The formula for the structure constants \ref{C1} in this case is:
\begin{eqnarray}
C^{\gamma}_{\alpha , \beta} = \sum_{m\in I_{g,g}(\alpha , \beta, \gamma)} \exp \lbrack   \pi \imath \frac{- \tau m^{2}}{ c^{3}(a+d)} \rbrack
\end{eqnarray}
The index set of the series is nonempty  only when 
$\alpha  \equiv d(\gamma-\beta) \mod c$:
\begin{eqnarray*}
I_{g,g}(\alpha , \beta, \gamma) \neq \emptyset &  \iff  & -c\gamma + c(a+d)\alpha  \equiv c(d^2 +bc)\gamma - c(a+d)d \beta \mod  c^2(a+d) \\
&  \iff  & -\gamma + (a+d)\alpha  \equiv (d^2 +bc)\gamma - (a+d)d \beta  \mod   c(a+d) \\
&  \iff  & (a+d)\alpha  \equiv (d^2 +bc +1 )\gamma - (a+d)d \beta  \mod    c(a+d) \\
&  \iff  & (a+d)\alpha  \equiv (d^2 + da )\gamma - (a+d)d \beta  \mod   c(a+d) \\
&  \iff  & (a+d)\alpha  \equiv d(d + a)\gamma - (a+d)d \beta  \mod    c(a+d) \\
&  \iff  &  \alpha  \equiv d(\gamma - \beta)  \mod   c \\
\end{eqnarray*} 
Thus $C^{\gamma}_{\alpha , \beta} = 0 $ if $ \alpha  \not\equiv d(\gamma - \beta) \mod c$.

Conversely if $ \alpha  \equiv d(\gamma - \beta) \mod c$ then 
\begin{eqnarray*}
I_{g,g}(\alpha , \beta, \gamma) 
&=& \{  n \in  \mathbb{Z} | n \equiv -c\gamma + c(a+d)\alpha    \mod c^2(a+d)    \} \\
&=& \{  n \in  \mathbb{Z} | n = -c\gamma + c(a+d)\alpha  +  m c^2(a+d) \text{ for some }   m \in  \mathbb{Z}  \}
\end{eqnarray*}

And thus 
\begin{eqnarray*}
C^{\gamma}_{\alpha , \beta} 
&  = &  \sum_{n \in  \mathbb{Z} } \exp [- \pi \imath \tau \frac{(-c\gamma + c(a+d)\alpha  +  n c^2(a+d) )^{2}}{ c^{3}(a+d)} ]\\
&  = &  \sum_{n \in  \mathbb{Z} } \exp [- \pi \imath \tau ( \frac{(-c\gamma + c(a+d)\alpha)^{2}}{ c^{3}(a+d)}+\frac{2n (-c\gamma + c(a+d)\alpha)}{c} + c(a+d)n^{2} )]\\
&  = &  \exp[- \pi \imath \tau \frac{( (a+d)\alpha  -\gamma)^{2}}{ c(a+d)}] \sum_{n \in  \mathbb{Z} } \exp [- 2 \pi \imath \tau ( (a+d)\alpha  -\gamma) n -  \pi \imath \tau c(a+d) n^2 ]\\
\end{eqnarray*}
The last series is a theta series \ref{thetaseries}:
\begin{eqnarray*}
C^{\gamma}_{\alpha , \beta}
&  = &  \exp[- \pi \imath \tau \frac{[ (a+d)\alpha  -\gamma]^{2}}{ c(a+d)}] \vartheta (- \tau ( (a+d)\alpha  -\gamma) ,- \tau c(a+d) ) \\  
\end{eqnarray*}
And we can write it as a theta constant with rational coefficients by 
taking, $\tau'= - \tau c(a+d)$ and $l= c(a+d)$. 
\begin{eqnarray*}
C^{\gamma}_{\alpha , \beta} 
&  = &  \exp[- \pi \imath \tau c(a+d) [\frac{ (a+d)\alpha  -\gamma}{ c(a+d)}]^2]  \vartheta ( -\tau c(a+d) \frac{ (a+d)\alpha  -\gamma}{ c(a+d)} ,- \tau
c(a+d) ) \\
&  = &  \exp[\pi \imath \tau' [\frac{ (a+d)\alpha  -\gamma}{ l}]^2]  \vartheta ( \tau '\frac{ (a+d)\alpha  -\gamma}{l},
\tau ') \\
&  = &  \vartheta_{\frac{(a+d)\alpha  -\gamma}{l}} ( \tau ').
\end{eqnarray*}

Now, by \ref{zerostheta} for $r,s\in \frac{1}{l}\mathbb{Z}$ the zeroes of $\vartheta_{r,s} ( z , \tau ') $ occur at the points of the form $(r+p + \frac{1}{2})\tau ' + (s + q + \frac{1}{2} )$ for $p,q \in \mathbb{Z}$. In particular the zeroes of $\vartheta_{r,0}$ are at points $(r+p + \frac{1}{2})\tau ' + (q + \frac{1}{2})$. Thus $\vartheta_{r} (\tau ') = \vartheta_{r,0} (0, \tau')\neq  0$ for all $r\in \frac{1}{l}\mathbb{Z}$ which proves the lemma. 
\end{proof}

The expression for the nonzero values of the structure constants in 
the proof of Lemma~\ref{lemmaconstantes} is crucial in all that follows. We state it as 
a corollary.
\begin{corollary} The nonzero values of the structure constants $C^{\gamma}_{\alpha , \beta}$ are
theta constants with rational characteristics depending on $g$. If $ \alpha  \equiv d(\gamma - \beta) \mod c$ then
\begin{eqnarray*}
C^{\gamma}_{\alpha , \beta} &  = &  \vartheta_{\frac{(a+d)\alpha  -\gamma}{l}} ( \tau ') \\
&  = & \vartheta_{\frac{(a+d)d(\gamma - \beta)  -\gamma}{l}} ( \tau ')  
\end{eqnarray*}
where $\tau'= - \tau c(a+d)$ and $l= c(a+d)$.
\end{corollary}
\begin{proof}
The first equality follows from the proof of  Lemma~\ref{lemmaconstantes}. For the second one note that the theta constant $\vartheta_{r} ( \tau ')$ only depend on the class of $r$ in $\frac{1}{l} \mathbb{Z} / \mathbb{Z}$ so we can replace $\alpha$ by $d(\gamma - \beta) \mod c$ in the last formula for the nonzero structure constants. 
\end{proof}

We can use Lemma~\ref{lemmaconstantes} to write the linear system $M$ corresponding to $m$ as a sum of $c$ independent systems. We do this by grouping the nonzero structure constants. Since $c$ and $d$ are relatively prime it follows that for $\gamma$ in a fixed congruence class $\mod c$ the value of $\beta $ determines the unique value of $\alpha$ for which  $C^{\gamma}_{\alpha , \beta} \neq 0$. 
\begin{notation}
\label{alpha}
Given $\mu, \beta \in \{ 1,...,c\}$ we denote by $\alpha(\mu,\beta)$ the unique representative $\mod c$ of $d(\mu-\beta)$ laying in $\in \{ 1,...,c\}$.
\end{notation}
\begin{lemma}
Let $M:\mathbb{C}^{c^2} \rightarrow \mathbb{C}^{c (a+d)}$ be given by $C^{\gamma}_{\alpha , \beta}$. Then $M$ is equivalent to $c$ independent systems $M(\mu):\mathbb{C}^{c} \rightarrow \mathbb{C}^{(a+d)}$, $\mu =1,...,c$, each one of rank $c-(a+d)$.   
\end{lemma}
\begin{proof}
Fix $\alpha$ and $\beta$ in $\{ 1,2, ...c \} $. Given $\gamma ,\gamma' \in \{ 1,2, ...c(a+d) \}$ we have by Lemma~\ref{lemmaconstantes} 
\begin{eqnarray*}
C^{\gamma}_{\alpha , \beta} \neq 0 \text{ and } C^{\gamma'}_{\alpha , \beta} \neq 0  
& \iff &  \alpha + d \beta \equiv d\gamma  \text{ and } \alpha + d \beta \equiv d\gamma'  \\
& \iff & d\gamma  \equiv d\gamma'  \mod c\\ 
& \iff & \gamma  \equiv \gamma' \mod c \\
\end{eqnarray*}

Therefore nonzero values of $C^{\gamma}_{\alpha , \beta}$  occur for values of $\gamma$ in the same congruence class $\mod c$ and we can arrange the system as $c$ independent systems of dimension $c^{2}\times (a+d)$, each one corresponding to the values of of $\gamma$ in the same congruence class $\mod c$. Again by Lemma~\ref{lemmaconstantes} each one of these systems will have only $c$ nontrivial columns corresponding to the values of $\alpha$ and $\beta$ satisfying 
$\alpha  \equiv d(\gamma - \beta) \mod c$. Leaving aside the zero structure constants we are left with the $c$ independent systems:  
\begin{eqnarray}
\label{Mmu}
M(\mu) = \left(
\begin{array}{cccc}
  C^{\mu}_{\alpha(\mu,1) , 1} & C^{\mu}_{\alpha(\mu,2) ,2} &...&C^{\mu}_{\alpha(\mu,c) , c} \\
  C^{\mu +c}_{\alpha(\mu,1) , 1} & C^{\mu+c}_{\alpha(\mu,2) ,2} &...&C^{\mu+c}_{\alpha(\mu,c) , c} \\
  . &  &  & \\
  . &  &  & \\
  C^{\mu +(a+d-1)c}_{\alpha(\mu,1) , 1} & C^{\mu+(a+d-1)c}_{\alpha(\mu,2) ,2} &...
  &C^{\mu +(a+d-1)c}_{\alpha(\mu,c) , c} \\
\end{array}
\right)
\end{eqnarray}
where $\mu \in \{ 1, 2, ...,c \}$.
By Theorem~\ref{poli2} $B_{g}(\theta,\tau)$  is generated by its degree 1 part (remember we assumed $c \geq a+d +2$). This in particular means that $\mathcal{H}_{g^2}$, its degree two part, is generated by products of elements in $\mathcal{H}_g$ so the multiplication map $m$ is surjective. At the level of the representing matrices this just means that $M$ has maximal rank. Also, since $M$ is the direct sum of the $M(\mu)$ each one of these must have maximal rank. 
\end{proof}
Componentwise we have  $M(\mu)_{i,j }= C^{\mu + (i -1)c}_{\alpha(\mu,j) , j}$. 
By Lemma~\ref{lemmaconstantes} the denominator in the characteristic of the corresponding 
theta constant is $(a+d)d(\mu -j) - (\mu + (i-1)c)$ so dividing out by $l=c(a+d)$ one gets $\frac{d\mu}{c}-\frac{dj}{c}-\frac{\mu}{l}-\frac{i}{a+d}+\frac{1}{a+d}$
i.e.
\begin{eqnarray}
\label{Mmu2}
M(\mu )_{i,j}= \vartheta_{q(\mu)-\frac{dj}{c}-\frac{i}{a+d}} ( \tau ')
\end{eqnarray}
Where $q(\mu)= \frac{d\mu}{c} -\frac{\mu}{l}+\frac{1}{a+d}$ is the term not depending on $i$ or $j$.
It is important to note that the characteristics giving the theta constants which appear as the coefficients of the $M(\mu)$ are all the same up to a shift. We can arrange these characteristics in a matrix $\Lambda=\Lambda(g)\in \mathcal{M}_{c\times(a+d)}(\frac{1}{l}\mathbb{Z}/ \mathbb{Z})$ given by $\Lambda_{i,j}= -\frac{dj}{c}-\frac{i}{(a+d)}$. The matrices $M(\mu)$ are then functions of $\tau$ and $\mu$ determined by $g$: 
\begin{eqnarray}
\label{Mmu3}
(\tau,\mu) \mapsto M(\mu)_{i,j}= \vartheta_{q(\mu)+\Lambda_{i,j}}( \tau') 
\end{eqnarray}

Each one of the kernels of the matrices $M(\mu)$ gives us a set of relations for $B_{g}(\theta,\tau)$. By the above discussion we see that these sets are independent. We write down a basis for the kernel of each $M(\mu)$ in terms of its minor determinants. First we introduce some notation:
\begin{notation}{}
Let $n >m$ and let $L:\mathbb{C}^{n} \rightarrow \mathbb{C}^{m}$ be a surjective linear map. 
Denote also by $L \in \mathcal{M}_{m\times n}(\mathbb{C})$ its matrix representation in the 
canonical basis. Given $i_1,i_2,...,i_m \in \{ 1,2,...,n\}$, $i_1<i_2<...<i_m$ we write
$| L| _{i_1,i_2,...,i_m}$ for the minor determinant corresponding to the columns 
$ i_1,i_2,...,i_m $.
\end{notation}

\begin{lemma}
\label{lemmalineal}
Let $n >m$ and let $L \in \mathcal{M}_{n\times m}(\mathbb{C})$. Assume the first $m$ columns of 
$L$ are linearly independent. For each $k=1,...,n-m$ let $v^{k}\in \mathbb{C}^{n}$ be defined by:
\begin{eqnarray}
\label{v0}
v^{k}_j  &= 
\begin{cases} 
| L| _{1,2,...,j-1,m+k,j+1,...,n} & if\,  1\leq j \leq m \\
-| L| _{1,2,...,m} & if \, j=m+k \\
0 & \text{otherwise}
\end{cases}
\end{eqnarray}
Then $\{ v^{k} | k=1,...,n-m\}$ forms a basis for $Ker(L)$
\end{lemma}
\begin{proof} 
Let $L$ be as above. Denote by $L^1,...,L^n \in \mathbb{C}^{m}$ its columns. Denote by $\tilde{L}\in \mathcal{M}_{m\times m}$ the  matrix corresponding to the first $m$ columns of $L$. 
For $k\in\{1,...,n-m\}$ let  $\tilde{X^k} \in \mathbb{C}^{m}$ be a solution of $\tilde{L} Y = -L^{m+k}$ we can complete $\tilde{X^k}$ to a vector $X^k \in \mathbb{C}^{n}$ by taking the remaining coordinates to be $0$ except for the $n+k$ coordinate which we set to $1$. Then $L X^k=0$ for all $k\in\{1,...,n-m\}$. Also, it is clear from the construction that the $X^k$ are linearly independent. In this way we may construct a basis for the
kernel of $L$. Now, we solve each one of the systems $\tilde{L} Y = -L^{m+k}$ using  Cramer's rule. After clearing denominators in the solution and completing to a vector in 
$\mathbb{C}^n$ we get the vectors $v^{k}$ in \ref{v0}. 
\end{proof}

It will be convenient to view the minor determinants of $M(\mu)$ as functions of $\tau$:
\begin{definition}
Let $M(\mu)$ be given as in \ref{Mmu} let $i_1,...,i_{a+d} \in \{ 1,2,...,c\}$ with $i_1<i_2<...<i_{a+d}$. Define
\begin{eqnarray*}
F^{g, \mu}_{i_1,i_2,...,i_{a+d}}(\tau) 
&  = & |M(\mu)|_{i_1,i_2,...,i_{a+d}}  \\
& = & \sum_{\sigma \in S_{a+d}} \text{sgn} (\sigma) \prod_{k=1}^{a+d} M(\mu)_{\sigma (k),i_k} \\
& = &  \sum_{\sigma \in S_{a+d}} \text{sgn} (\sigma) \prod_{k=1}^{a+d}
\vartheta_{q(\mu)-\frac{d i_{k}}{c}-\frac{\sigma (k)}{(a+d)}} (  \tau ') \\
\end{eqnarray*}
\end{definition}

We will show later that for each $g$ and $\mu$ these are modular functions on $\tau$. By applying Lemma~\ref{lemmalineal} we can write now a explicit presentation of  $B_g(\tau,\theta)$:
\begin{theorem}
\label{precentacion1}
Given $\mu \in \{ 1, 2, ...,c \}$ and  $k \in \{ 1,..., c-a-d \}$ let $v^{\mu,k} = v^{\mu,k}(\tau)\in \mathbb{C}^{c}$ be given by 

\begin{eqnarray}
\label{v1}
v^{\mu,k}_j = v^{\mu,k}_j(\tau) =
\begin{cases} 
 F^{g, \mu}_{1,2,...,j-1,a+d+k,j+1,...,a+d}(\tau) & if \,  1\leq j \leq a+d \\
 F^{g, \mu}_{1,2,...,a+d} (\tau) & if\,  j=m+k \\
 0 &\text{otherwise}
\end{cases}
\end{eqnarray}
Then the algebra $B_g(\tau,\theta)$ is generated by elements $x_1,...,x_c$ of degree $1$ subject to relations $\mathit{f}^{\mu}_{k}=0$ where:
\begin{eqnarray}
\mathit{f}^{\mu}_{k}= v^{\mu,k}_1 x_{\alpha(\mu,1)} x_1 +...+v^{\mu,k}_c x_{\alpha(\mu,c)} x_c
\end{eqnarray}
\end{theorem}
\begin{proof}
Each one of the matrices $M(\mu)$ has maximal rank equal to $a+d$. Therefore there are $a+d$ linearly independent columns and we can reorder them in order to apply Lemma~\ref{lemmalineal}.
We let $x_1,...,x_c$ be the generators of the tensor algebra $\mathcal{T}$ of $\mathcal{H}_g$  corresponding to the basis $\varphi_1,...,\varphi_c$. Thus $\mathcal{T}=\mathbb{C}\langle x_1,...,x_c \rangle$. For each $\mu \in \{ 1, 2, ...,c \}$ Lemma~\ref{lemmalineal} gives us a basis for the kernel of $M(\mu)$ which corresponds by Lemma~\ref{lemmafacil} to a set of defining relations for $B_g(\tau,\theta)$. 
\end{proof}

\begin{example}
Let  
$$g= \left(
\begin{array}{cc}
4 & -1 \\
5 & -1 \\
\end{array}
\right)$$ 
The eigenvalues of $g$ are 
$$
\lambda^{+} = \frac{3 - {\sqrt{5}}}{2}, \quad  \lambda^{-} = \frac{3 + {\sqrt{5}}}{2}
$$ 
and the fixed points of $g$ are 
$$
\theta = \frac{5 - {\sqrt{5}}}{10}, \quad  \theta' = \frac{5 + {\sqrt{5}}}{10}.
$$
Fix now a complex structure $\tau$ on $\mathcal{A}_{\theta}$ and consider the corresponding connection $ \bar{\nabla}_0$ on the $\mathcal{A}_{\theta}$-bimodule $E_{-1,5}(\theta)$.
Since 
$$g^{2} = \left(
\begin{array}{cc}
  11 & -3 \\
  15 & -4 \\
\end{array}
\right)$$ 
we have that $\mathcal{H}_g \simeq \mathbb{C}^{5}$ and $\mathcal{H}_{g^2}\simeq \mathbb{C}^{15}$.
After choosing a basis the multiplication map $m:\mathcal{H}_{g}\otimes \mathcal{H}_{g}\rightarrow \mathcal{H}_{g^2}$ is represented by a matrix $M\in \mathcal{M}_{15,25}(\mathbb{C})$. We take as above $\{ \varphi_{\alpha} \otimes \varphi_{\beta}| \alpha, \beta = 1,...,5 \}$ as basis for $\mathcal{H}_{g}\otimes \mathcal{H}_{g}$ and $\{ \psi_{\gamma}| \gamma = 1, ..., 15 \}$ as basis for $\mathcal{H}_{g^2}$ so that $M$ is the matrix corresponding to the structure constants 
$$
C^{\gamma}_{\alpha , \beta} \quad =\quad 
\begin{cases} 
 \quad \vartheta_{\frac{3 \beta -4 \gamma }{15}} (-15 \tau ) 
 \quad if\quad  \alpha  \equiv d(\gamma - \beta) \mod c \\
 \quad 0 \quad \text{otherwise}
\end{cases}
$$
We write it as $M\simeq M(1)\oplus  M(2)\oplus  M(3)\oplus  M(4)\oplus  M(5)$ where the elements of $M(\mu )\in  \mathcal{M}_{3,5}(\mathbb{C})$ are given by 
$$
M(\mu )_{i,j}= \vartheta_{\frac{5 - 4 \mu}{15}+\frac{j}{5}-\frac{i}{3}} ( - 15 \tau )
$$
The corresponding matrix of characteristics is:
$$
\Lambda = \frac{1}{15}\left(
\begin{matrix}
2 & 14 & 11 & 8 & 5 \cr 7 & 4 & 1 & 13 & 10 \cr 12 & 9 & 6 & 3 & 0 \cr  
\end{matrix}
\right)
$$
Each $\mu\in\{1, ...,5 \}$ gives us a set of $2=c-a-d$ relations corresponding to a basis for the kernel of the matrix $M(\mu)$. Thus in this case $B_{g}(\theta,\tau)$ is a quadratic algebra with $5$ generators of degree $1$ and ten quadratic relations. The corresponding functions of $\tau$ determining the coefficients of the relations are the minors of each $M(\mu)$. In this case 
all minors are non vanishing. Each ordered triple $i_1,i_2,i_3 \in \{ 1,2,3,4,5\}$ determines a function of $\tau$ corresponding to the determinant of the matrix formed by the three corresponding  columns. 
\begin{eqnarray*}
F^{g, \mu}_{i_1,i_2,i_{3}}(\tau) 
&  = & |M(\mu)|_{i_1,i_2,i_{3}}  \\
& = &  \sum_{\sigma \in S_{3}} \text{sgn} (\sigma) \prod_{k=1}^{3}
\vartheta_{\frac{5 - 4 \mu}{15}-\frac{d i_{k}}{5}-\frac{\sigma (k)}{3}} ( -15 \tau ) \\
\end{eqnarray*}
Thus each one of the coefficients of the defining relations is a sum of triple products of theta constants. For example, taking $\mu =1$ we get the two relations:  
$$
\mathit{f}^{1}_{1}= 
F^{g, 1}_{4,2,3}(\tau) x_{5} x_1 
+ F^{g, 1}_{1,4,3}(\tau) x_{1} x_2 
+F^{g, 1}_{1,2,4}(\tau) x_{2} x_3 
-F^{g, 1}_{1,2,3}(\tau) x_{3} x_4 
$$
and
$$
\mathit{f}^{1}_{2}= 
F^{g, 1}_{5,2,3}(\tau) x_{4} x_1 
+ F^{g, 1}_{1,5,3}(\tau) x_{5} x_2 
+F^{g, 1}_{1,2,5}(\tau) x_{1} x_3 
-F^{g, 1}_{1,2,3}(\tau) x_{3} x_5 .
$$ 
\end{example}


\section{Rationality of $B_g(\tau,\theta)$}
\label{sec4}

On this section and the next one we use the presentation of $B_g(\tau,\theta)$ in terms of theta functions given in Theorem~\ref{precentacion1} in order to derive some results about the rationality of $B_g(\tau,\theta)$ and its behavior for different values of $\tau$. In each case the main idea is to rescale the relations by the appropriate factors. As we saw above the relations defining $B_g(\tau,\theta)$ are given by a linearly independent generating set for the kernel of the multiplication map. Multiplying each element of this basis by a nonzero constant we still obtain a basis. Therefore we are free to multiply each one of the 
relations in Theorem~\ref{precentacion1} by a nonzero constant and still get a set of defining relations for $B_g(\tau,\theta)$.

Let $E_{\tau'}$ be the elliptic curve $\mathbb{C}/ (\mathbb{Z}\oplus\tau' \mathbb{Z})$. Let $j(\tau')$ be the absolute invariant of $E_{\tau'}$ and denote by $k'=\mathbb{Q}(j(\tau'))$ the minimal field of definition of $E_{\tau'}$. Since the group structure on the elliptic curve $E_{\tau'}$ is given in terms of algebraic maps any torsion point on $E_{\tau'}$ is defined over a finite algebraic extension of the field $k'$. If we denote by $E_{\tau'}[N]$ the set of $N$-torsion points of $E_{\tau'}$ and let $k'_{N}$ be the field extension of $k'$ over which they are defined we obtain a finite algebraic extension $k' \hookrightarrow  k'_{N}$ (in fact it is a Galois extension c.f. \cite{lang} Chapter2 \S 1).
Since for any $r\in \mathbb{Q}$ we have that $\vartheta_{r}(\tau') = \vartheta_{r,0}^{\alpha}(\tau')$ we can apply Theorem~\ref{teotata3} to study the rationality of the coefficients of the defining relations of $B_g(\tau,\theta)$. Let $l=c(a+d)$, given  $r \in \frac{1}{l}\mathbb{Z}$ the point $x(r,0)_{2 l}$ is a torsion point of order dividing $2l^2$ and is therefore rational over $k'_{2l^2}$. Since $E_{\tau'}$ and its basic line bundle $\mathcal{L}$ are both defined over $k'_{2l^2}$ we get: 
\begin{theorem}
\label{values1}
For all $r \in \frac{1}{l}\mathbb{Z}$ 
\begin{eqnarray}
\frac{\vartheta_{r}(\tau') }{\vartheta(\tau')} \in k'_{2l^2}
\end{eqnarray}
\end{theorem}

As in the previous section let $\tau'= -l \tau$, denote by $E_{\tau}$ the corresponding elliptic curve $\mathbb{C}/ (\mathbb{Z}\oplus\tau \mathbb{Z})$. 
The minimal field of definition of the elliptic curve $E_{\tau}$ is 
$k=\mathbb{Q}(j(\tau))$. Since $k'_{2l^2}$ is a finite algebraic extension of $k'$ Theorem~\ref{values1} implies that 
$\frac{\vartheta_{r}(\tau') }{\vartheta(\tau') }$ is algebraic over $k'$.
We want to show that the values 
$\frac{\vartheta_{r}(\tau') }{\vartheta(\tau') }$ are also algebraic over $k$. For this we must study the relation between $k$ and $k'_{2l^2}$. The main point here is that the absolute invariants $j(\tau)$ and $j(\tau')$ are related by an algebraic equation with integer coefficients:
\begin{theorem}\emph{(c.f \cite{lang} Chapter 5 \S 3)}
\label{teoinv}
Let $\Phi_{n}(T,T')\in \mathbb{Z}[T,T']$ be the modular polynomial of order $n$ and let $E_{1}, E_{2}$ be two elliptic curves over $\mathbb{C}$. Then $\Phi_{n}(j(E_{1}),j(E_{2}))=0$ if and only if  there exist an isogeny 
$E_{1} \rightarrow E_{2}$ with cyclic kernel of degree $n$. 
\end{theorem}

\begin{corollary}
\label{postlang}
Let $l>1$ be an integer and let $r\in \frac{1}{l}\left( \mathbb{Z} \right)$. Then 
\begin{eqnarray*}
\frac{\vartheta_{r}(\tau') }{\vartheta(\tau') }= \frac{\vartheta_{r}(- l\tau)}{\vartheta(- l \tau) }
\end{eqnarray*} 
is algebraic over $k$.
\end{corollary}
\begin{proof}
Let $E_{\tau}$ and $E_{\tau'}$ be as above. Multiplication by $l$ induces an isogeny $E_{\tau} \rightarrow E_{\tau'}$. The kernel of this isogeny is $(\mathbb{Z}\oplus\tau \mathbb{Z})/(\mathbb{Z}\oplus l \tau \mathbb{Z}) \simeq \mathbb{Z}/ l \mathbb{Z}$. From Theorem~\ref{teoinv} it follows that $j(\tau')$ is algebraic over $k=\mathbb{Q}(j(\tau))$. Taking into account that $k'_{2l^2}$ is algebraic over $\mathbb{Q}(j(\tau'))$ we get a tower of algebraic extensions:
\begin{eqnarray*}
\mathbb{Q}(j(\tau)) \hookrightarrow \mathbb{Q}(j(\tau'),j(\tau)) \hookrightarrow k'_{2l^2}(j(\tau)).
\end{eqnarray*} 
The theorem then follows from Theorem~\ref{values1}
\end{proof}

We can state now the main theorem of this section:
\begin{theorem}
\label{superteorema}
Let $E_{\tau}$ be the elliptic curve $\mathbb{C}/ (\mathbb{Z} \oplus\tau \mathbb{Z})$. Let $k$ be its minimal field of definition. Then the algebra $B_{g}(\theta,\tau)$ admits a rational presentation over a finite algebraic extension of $k$.
\end{theorem}
\begin{proof}
Let $\mathbb{K}= k'_{2l^2}(j(\tau))$. By Corollary~\ref{postlang}, $\mathbb{K}$ is an algebraic extension of $k$ and $\frac{\vartheta_{r}(\tau') }{\vartheta(\tau') }$ belongs to $\mathbb{K}$ for every  $r\in \frac{1}{l}\mathbb{Z}$.

Multiplying each one of the relations in Theorem~\ref{precentacion1} by
$\left[\vartheta(\tau')\right]^{-(a+d)}$ we get a new basis $\{\tilde{\mathit{f}}^{\mu}_{k}\}$ for the ideal $\mathcal{R}$ giving the relations of the $B_g(\tau,\theta)$:
\begin{eqnarray*}
\tilde{\mathit{f}}^{\mu}_{k}= \tilde{v}^{\mu,k}_1 x_{\alpha(\mu,1)} x_1 +...
+\tilde{v}^{\mu,k}_c x_{\alpha(\mu,c)} x_c
\end{eqnarray*}
where $\tilde{v}^{\mu,k}_j =\left[\vartheta(\tau')\right]^{-(a+d)}v^{\mu,k}_c$. Each 
$\tilde{v}^{\mu,k}_j $  has the form:
$$
\left[\vartheta( \tau')\right]^{-(a+d)} F^{g, \mu}_{i_1,...,i_{a+d}}(\tau)  =
\left[\vartheta( \tau')\right]^{-(a+d)}  \sum_{\sigma \in S_{a+d}} \text{sgn} 
(\sigma) \prod_{k=1}^{a+d}
\vartheta_{q(\mu)-\frac{d i_{k}}{c}-\frac{\sigma (k)}{(a+d)}} (\tau ') 
$$ for some $i_1,...,i_{a+d}\in \{ 1,...,a+d \}$. Therefore $\tilde{v}^{\mu,k}_j$ belongs to $\mathbb{K}$ being a rational combination of values of the form  $\frac{\vartheta_{r}(\tau') }{\vartheta(\tau') }$ .
\end{proof} 

The field $\mathbb{K}$ is generated over $\mathbb{Q}$ by $j(\tau),j(\tau')$ and the values of the Weierstrass functions $\wp$ and $\wp'$ at the points $\frac{p}{2l^2} \tau ' + \frac{q}{2l^2}$ with $p,q\in\mathbb{Z}$. From the preceding discussion that $B_{g}(\theta,\tau)$ is defined over $\mathbb{K}$.

In particular we can restrict our field of scalars to $\mathbb{K}$ to obtain a noncommutative $\mathbb{K}$-algebra:
\begin{eqnarray}
\label{LAdefinicion2}
B_g(\tau,\theta)_{\mathbb{K}}=\mathbb{K} \langle x_1, ..., x_c \rangle / \mathcal{R}
\end{eqnarray}

An important question is whether the Galois group $Gal(\mathbb{K}/k)$ leaves $\mathcal{R}$ invariant and being this the case how does it act on $B_g(\tau,\theta)_{\mathbb{K}}$.


\section{Special values of $\tau$}

Starting with  Theorem~\ref{superteorema} we can study the properties 
of $B_g(\tau,\theta)$ for special values of $\tau$ giving interesting fields of definition. 
We use the same notation as in Section~\ref{sec4}.
\subsection{Subfields of $\mathbb{R}$}
We start this section with a simple, yet very useful, remark:
\begin{proposition}
\label{real}
Let $\tau\in \imath \mathbb{R}$. Then $B_g(\tau,\theta)$ is defined over a subfield of $\mathbb{R}$. 
\end{proposition}
\begin{proof}
$B_g(\tau,\theta)$ is defined over the field generated over $\mathbb{Q}$ by the values of the 
theta constants $\vartheta_{r_{i}}(\tau')$. By the series expression \ref{thetaseries} we know 
that the theta constant $\vartheta_{r_{i}}$ takes real values on $\imath \mathbb{R}$. Since $\tau\in \imath \mathbb{R}$ implies $\tau' \in \imath \mathbb{R}$ the proposition follows.
\end{proof}    

\subsection{Number fields}
The most interesting family of examples comes from elliptic curves defined over number fields. In the case $j(\tau)$ is algebraic over $\mathbb{Q}$ the field $\mathbb{K}= k'_{2l}(j(\tau))$ is a number field. As before we consider the algebra 
$$
B_g(\tau,\theta)_{\mathbb{K}} =  \mathbb{K}\langle x_1,...,x_c \rangle / \mathcal{R}
$$
obtained by restriction of scalars from $\mathbb{C}$ to $\mathbb{K}$. Let now $\mathcal{O}_{\mathbb{K}}$ be the ring of integers of $\mathbb{K}$. Since $\mathbb{K}$ is the field of fractions of $\mathcal{O}_{\mathbb{K}}$ we can clear denominators in each one of the defining relations of $B_g(\tau,\theta)_{\mathbb{K}}$ and obtain a basis $\{ \bar{\mathit{f}}^{\mu}_{k} \}$ of $\mathcal{R}$ of the form
$$
\bar{\mathit{f}}^{\mu}_{k} = \bar{v}^{\mu,k}_1 x_{\alpha(\mu,1)} x_1 +...+\bar{v}^{\mu,k}_c x_{\alpha(\mu,c)} x_c
$$
with $\bar{v}^{\mu,k}_j$ in $\mathcal{O}_{\mathbb{K}}$. In particular we can consider the various 
reductions corresponding to different finite places $\mathfrak{P}$ of $\mathbb{K}$. If the coefficient $\bar{v}^{\mu,k}_j$ are nonzero modulo $\mathfrak{P}$ it makes sense to talk about the $\mathcal{O}_{\mathbb{K}} /\mathfrak{P}$ algebra given by the relations $\bar{\mathit{f}}^{\mu}_{k} \mod \mathfrak{P}$.

Among number fields some cases deserve particular attention. The first case comes from taking $\tau = \sqrt{-D}$ where $D$ is a positive integer. Since the $j$ invariant of an elliptic curve with complex multiplication is algebraic we get from Proposition~\ref{real}: 
\begin{corollary}
Let $\tau = \sqrt{-D}$ be a generator in $\imath \mathbb{R}$ of the quadratic imaginary 
field $\mathbb{Q}(\sqrt{-D})$. Then $B_g(\tau,\theta)$ is defined over a real algebraic extension of $\mathbb{Q}$. 
\end{corollary}

In particular we can choose $\tau = \sqrt{-D}$ where $D$ is a positive integer such that $\mathbb{Q}(\theta)= \mathbb{Q}(\sqrt{D})$.

Finally, another case that could give rise to interesting structures comes from considering an elliptic curve defined over $\mathbb{Q}(\theta)$. For instance take $\tau$ with absolute invariant $j(\tau)=\theta$. Then one gets defining relations rational over $\mathbb{Q}(\theta)$ which is the field we ultimately want to study.


\section{Linear basis for $B_g(\tau,\theta)$}

Let $\mathbb{K}$ be a subfield of $\mathbb{C}$ over which $B_g(\tau,\theta)$ is defined. We would like to construct interesting linear functionals on $B_g(\tau,\theta)_{\mathbb{K}}= \mathbb{K} \langle x_1, ..., x_c \rangle / \mathcal{R}$ with values on $\mathbb{K}$. Provided we are given a linear basis for $B_g(\tau,\theta)_{\mathbb{K}}$ over $\mathbb{K}$ we can  define linear functionals by prescribing the values of the elements in such a basis. The aim of this section is to study the natural linear basis corresponding to the presentations of $B_g(\tau,\theta)$ in terms of generators and relations given in Theorem~\ref{superteorema}. For this purpose the theory of noncommutative Gr\"obner basis for two sided ideals on the free algebra $\mathbb{K} \langle x_1, ..., x_c \rangle$ provides the right framework. The general idea is an extrapolation of Gaussian reduction to infinite dimensions (c.f. \cite{Mora}, see also 
\cite{Polibook}).

We view $\mathbb{K}\langle x_1, ..., x_c \rangle$ as the semigroup algebra corresponding to the free semigroup $\mathbf{S}$  generated by $\{x_1, ..., x_c \}$. The semigroup $\mathbf{S}$ is graded by total degree and it becomes an ordered semigroup by imposing on it the {\it deglex} order. Thus given $t_1,t_2 \in \mathbf{S}$ we say that $t_1<t_2$ if either $deg(t_1)<deg(t_2)$ or $deg(t_1)=deg(t_2)$ and there exist $l,r_1,r_2 \in S$ such that $t_1 =  l x_{j} r_1  $ and $t_2 =  l x_{i} r_2$ with $j<i$. Once this order is given every element $f$ of $\mathbb{K}\langle x_1, ..., x_c \rangle$ has a well defined leading term $T(f)\in \mathbf{S}$ and a leading coefficients $lc(f)\in \mathbb{K}$.

As in the previous section let 
$g= \left(
\begin{array}{cc}
  a & b \\
  c & d \\
\end{array}
\right)$ be a matrix satisfying satisfy \ref{condiciones1} and \ref{condiciones2}. For $\mu,j \in \{ 1, 2, ...,c \}$ and  $k \in \{ 1,..., c-a-d \}$ let $\tilde{v}^{\mu,k}_j$ be defined as in Theorem~\ref{superteorema}. In order to keep track of the ordering of the terms in the relations defining $B_g(\tau,\theta)$ we interchange the roles of the two factors in the degree two part of the free algebra $\mathbb{K}\langle x_1, ..., x_c \rangle$. Also, we divide each one of the relations $\{\tilde{\mathit{f}}^{\mu}_{k}\}$ by its leading coefficient. Accordingly we define a
basis $\{\hat{\mathit{f}}^{\mu}_{k}\}$ for $\mathcal{R}$, the ideal of relations of $B_g(\tau,\theta)$, by:
\begin{eqnarray*}
\hat{\mathit{f}}^{\mu}_{k} &=&  x_{a+d+k}  x_{\alpha(\mu,a+d+k)}+ 
\frac{\tilde{v}^{\mu,k}_{a+d}}{\tilde{v}^{\mu,k}_{a+d+k}} x_{a+d}  x_{\alpha(\mu,a+d)} +...+
\frac{\tilde{v}^{\mu,k}_{1}}{\tilde{v}^{\mu,k}_{a+d+k}} x_{1}  x_{\alpha(\mu,1)}.
\end{eqnarray*}
 
In the above expression the terms of the relations appear in decreasing order with $T(\hat{\mathit{f}}^{\mu}_{k})= x_{a+d+k}  x_{\alpha(\mu,a+d+k)}$ and $lc(\hat{\mathit{f}}^{\mu}_{k})=1$. Consider now the following decomposition of $\mathbf{S}_{2}$, the set of degree two elements in $\mathbf{S}$:
\begin{eqnarray*}
\mathbf{S}_{2}' &=& \{  x_{i_1} x_{i_2}  \in  \mathbf{S}_{2} \; | \; i_1 \leq a+d \} \\
\mathbf{S}_{2}'' &=& \{  x_{i_1} x_{i_2}  \in  \mathbf{S}_{2} \; | \; i_1 > a+d \}.
\end{eqnarray*}

Then $t'<t''$ for any $t' \in \mathbf{S}_{2}', t''\in \mathbf{S}_{2}''$ and the defining relations of $B_g(\tau,\theta)$ have the form 
\begin{eqnarray}
t_{k}^{\mu}  &=& \sum_{s=1}^{a+d} c_{s}^{\mu,k} t_{s}^{'\mu} \qquad \;  
t_{s}^{'\mu}  \in \mathbf{S}_{2}' \; , t_{k}^{\mu}  \in \mathbf{S}_{2}''
\end{eqnarray}
with  $t_{k}^{\mu}= x_{a+d+k}  x_{\alpha(\mu,a+d+k)}$, 
$t_{s}^{'\mu} =  x_{s}  x_{\alpha(\mu,s)}$ and $c_{s}^{\mu,k} = - \frac{\tilde{v}^{\mu,k}_{s}}{\tilde{v}^{\mu,k}_{a+d+k}}$. In particular $t_{1}^{'\mu} < t_{2}^{'\mu}<...< t_{a+d}^{'\mu} < t_{k}^{\mu}$.

Denote by $\mathbf{S}_{n}$ the set of elements of degree $n$ in $\mathbf{S}$ and set 
\begin{eqnarray}
\mathbf{S}_{n}' &=&\{  x_{i_1} x_{i_2} ... x_{i_n}  \in  \mathbf{S}_{n} \; |
 \; i_1,i_2,...,i_{n-1} \leq a+d \} 
\end{eqnarray}

The set  $\mathbf{S}_{n}' = \cup_n \mathbf{S}_{n}'$ linearly spans $B_g(\tau,\theta)$ since any element $ t  = x_{j_1} x_{j_2} ... x_{j_n}\in \mathbf{S}_{n} \setminus \mathbf{S}_{n}'$ can be expressed by smaller elements modulo $\mathcal{R}_{n}$. 

We want to extract a basis from the set $\mathbf{S}_{n}' = \cup_n \mathbf{S}_{n}'$. By cardinality conditions we can see that this set is redundant. The cardinality of $\mathbf{S}_{n}'$ is $(a+d)^{n-1}c$ while the $n$-th graded component of $B_g(\tau,\theta)$ has dimension $deg(g^{n})= dim \mathcal{H}_{g^{n}}$. The  Hilbert series for $B_g(\tau,\theta)$ is given by (c.f \cite{Po}):

\begin{eqnarray}
\label{Hilbertseries}
h_{B_g(\tau,\theta)}(t)= \frac{1 + \left(  c -a-d \right) \,t + t^2}
  {1 - \left( a + d \right) \,t + t^2}.
\end{eqnarray}

Starting with the linear generating set $\mathbf{S}'$ and 
the Hilbert series $h_{B_g(\tau,\theta)}(t)$ we can find a linear basis for each graded 
piece of $B_g(\tau,\theta)$ by extracting a minimal set of linearly dependent elements 
from each $\mathbf{S}_{n}'$. The redundant elements will correspond to leading terms of elements of the ideal $\mathcal{R}$. A basis for the semigroup ideal $T(\mathcal{R}) = \{ T(f) \;| \; f\in  \mathcal{R} \}$ is provided by Buchberger's algorithm. Below we state the main parts of the formalism in our context (c.f. \cite{Mora}).

In general, given a two sided ideal $\mathcal{I}$ of $\mathbb{K}\langle x_1, ..., x_c \rangle$ we can consider the semigroup ideal formed by its leading terms $T(\mathcal{I})\subset  \mathbf{S}$ and  its complement on $ \mathbf{S}$, $O(\mathcal{I}):=  \mathbf{S} \setminus T(\mathcal{I})$. Then we have:
\begin{itemize}
\item $\mathbb{K}\langle x_1, ..., x_c \rangle = \mathcal{I} \oplus Span_{\mathbb{K}}(O(\mathcal{I}))$
\item $\mathbb{K}\langle x_1, ..., x_c \rangle / \mathcal{I} \simeq Span_{\mathbb{K}}(O(\mathcal{I}))$
\end{itemize}

A generating set $G\subset \mathcal{I}$ such that the semigroup ideal $T(G)$ generated by $\{ T(g)|g\in G\}$ coincides with  $T(\mathcal{I})$ is called a Gr\"obner basis for $\mathcal{I}$. 

For the ideal $\mathcal{R}$ we have that $O(\mathcal{R})_2= \mathbf{S}_{2}'$ and $O(\mathcal{R})_n \subset \mathbf{S}_{n}'$. Given a Gr\"obner basis $G$ for $\mathcal{R}$ the set of redundant elements of  $\mathbf{S}_{n}'$ will be $T(G)_{n}\cap \mathbf{S}_{n}'$. Buchberger's algorithm provides a way to find a Gr\"obner basis $G$ for $\mathcal{R}$ starting with the set of generators $\{ \mathit{f}^{\mu}_{k} \}$.
\begin{remark}
Since $B_g(\tau,\theta)$ has exponential growth we should not expect $\mathcal{R}$ to have a finite Gr\"obner basis. Still an infinite Gr\"obner basis $G$ exist and Buchberger's algorithm provides a way to compute its elements of a given degree.  
\end{remark}

Buchberger's algorithm is a infinite dimensional analog of Gaussian reduction. First one reduces the problem of finding a Gr\"obner basis to a linear algebra problem. Then the computation of the corresponding linear basis is done by a combinatorial manipulation of the vectors that takes into account order of the generators in $\mathbf{S}$. 
Let $V$ be a linear subspace of $\mathbb{K}\langle x_1, ..., x_c \rangle$. A linearly generating set $B$ of $V$ such that $T(V)=T(B)$ is called a Gauss generating set. A linearly basis $B$ of $V$ such that $T(V)=T(B)$ is called a Gauss Basis. $G$ is a Gr\"obner basis of a two sided ideal $\mathcal{I}$ of $\mathbb{K}\langle x_1, ..., x_c \rangle$ if and only if $\mathcal{G}= \{ lfr  \; | \; l,r \in S; \, f\in G  \}$ is a Gauss basis for $\mathcal{I}$.

Assume we start with a ordered set of generators $G' = \{ f_1,...,f_s \}$ of the ideal $\mathcal{I}$. Then $\mathcal{G}'= \{ lfr  \; | \; l,r \in \mathbf{S}; \, f\in G  \}$ is an ordered generating set for $\mathcal{I}$. From $\mathcal{G}'$ we can extract a canonical Gauss basis for $\mathcal{I}$ by taking a set of linearly independent elements $\mathcal{G}$ such that: 
\begin{itemize}
\item $T(\mathcal{G}) = T(\mathcal{G}')$
\item If $v\in \mathcal{G}$, $w\in \mathcal{G}'$ are such that $T(v)=T(w)$ then $v\leq w$
\end{itemize}

One has to consider the fact that different elements in $\mathcal{G}'$ may have the same leading term. This self obstructions have to be taken into account inductively. 
Given  $f\in \mathbb{K}\langle x_1, ..., x_c \rangle$ we say that $h$ is a {\emph normal form of $f$ with respect to $G'$} if:
\begin{itemize}
\item $f-h \in \mathcal{I}$. 
\item Either $h = 0$ or $T(h) \notin  T(G')$.
\end{itemize}

Given $j=1,...,s$ for each pair $(l,r)\in \mathbf{S}\times \mathbf{S}$ such that  
$$
lT(g_j)r \in (T(g_1),T(g_2),...,T(g_{j-1}))
$$ or  $lT(g_j)= T(g_j) r$ the normal form of the element $lT(g_j)r$ with respect to $G'$ must be added to the set $G'$. Choosing at each stage a minimal iredundant set of pairs $(l,r)$ and adding the corresponding normal forms to $G'$ we get a Gr\"obner basis $G$ of $\mathcal{I}$.

When we apply this procedure to $\{\tilde{\mathit{f}}^{\mu}_{k}\}$ we can algorithmically compute for each degree $n$ a set of obstructions in $\{ lfr  \; | \; l,r \in \mathbf{S} ; \, f\in \mathcal{R} \; , deg(lfr) = n  \}$ together with their normal forms having degree $n$. The leading terms of this normal forms correspond to the elements of $\mathbf{S}_{n}'$ which are redundant.


\begin{example}
To illustrate the above ideas we will compute the linear basis of the degree three part of the ring considered in Example~\ref{superejemplo} below. We start with the following basis for the ideal of relations $\mathcal{R}$:
\begin{eqnarray*}
\hat{\mathit{f}}^{1}_{1} &=&  x_{5}  x_{1}+ \frac{1}{\tilde{v}^{1,1}_{5}} (
\tilde{v}^{1,1}_{4} x_{4}  x_{2} +
\tilde{v}^{1,1}_{3} x_{3}  x_{3} +
\tilde{v}^{1,1}_{2} x_{2}  x_{4} +
\tilde{v}^{1,1}_{1} x_{1}  x_{5}) 
\end{eqnarray*}
\begin{eqnarray*}
\hat{\mathit{f}}^{2}_{1} &=&  x_{5}  x_{2}+ \frac{1}{\tilde{v}^{2,1}_{5}} (
\tilde{v}^{2,1}_{4} x_{4}  x_{3} +
\tilde{v}^{2,1}_{3} x_{3}  x_{4} +
\tilde{v}^{2,1}_{2} x_{2}  x_{5} +
\tilde{v}^{2,1}_{1} x_{1}  x_{6}) 
\end{eqnarray*}
\begin{eqnarray*}
\hat{\mathit{f}}^{3}_{1} &=&  x_{5}  x_{3}+ \frac{1}{\tilde{v}^{3,1}_{5}} (
\tilde{v}^{3,1}_{4} x_{4}  x_{4} +
\tilde{v}^{3,1}_{3} x_{3}  x_{5} +
\tilde{v}^{3,1}_{2} x_{2}  x_{6} +
\tilde{v}^{3,1}_{1} x_{1}  x_{1}) 
\end{eqnarray*}
\begin{eqnarray*}
\hat{\mathit{f}}^{4}_{1} &=&  x_{5}  x_{4}+ \frac{1}{\tilde{v}^{4,1}_{5}} (
\tilde{v}^{4,1}_{4} x_{4}  x_{5} +
\tilde{v}^{4,1}_{3} x_{3}  x_{6} +
\tilde{v}^{4,1}_{2} x_{2}  x_{1} +
\tilde{v}^{4,1}_{1} x_{1}  x_{2}) 
\end{eqnarray*}
\begin{eqnarray*}
\hat{\mathit{f}}^{5}_{1} &=&  x_{5}  x_{5}+ \frac{1}{\tilde{v}^{5,1}_{5}} (
\tilde{v}^{5,1}_{4} x_{4}  x_{6} +
\tilde{v}^{5,1}_{3} x_{3}  x_{1} +
\tilde{v}^{5,1}_{2} x_{2}  x_{2} +
\tilde{v}^{5,1}_{1} x_{1}  x_{3}) 
\end{eqnarray*}
\begin{eqnarray*}
\hat{\mathit{f}}^{6}_{1} &=&  x_{5}  x_{6}+ \frac{1}{\tilde{v}^{6,1}_{5}} (
\tilde{v}^{6,1}_{4} x_{4}  x_{1} +
\tilde{v}^{6,1}_{3} x_{3}  x_{2} +
\tilde{v}^{6,1}_{2} x_{2}  x_{3} +
\tilde{v}^{6,1}_{1} x_{1}  x_{4})
\end{eqnarray*}
\begin{eqnarray*}
\hat{\mathit{f}}^{1}_{2} &=&  x_{6}  x_{1}+ \frac{1}{\tilde{v}^{1,2}_{5}} (
\tilde{v}^{1,2}_{4} x_{4}  x_{3} +
\tilde{v}^{1,2}_{3} x_{3}  x_{4} +
\tilde{v}^{1,2}_{2} x_{2}  x_{5} +
\tilde{v}^{1,2}_{1} x_{1}  x_{6}) 
\end{eqnarray*}
\begin{eqnarray*}
\hat{\mathit{f}}^{2}_{2} &=&  x_{6}  x_{2}+ \frac{1}{\tilde{v}^{2,2}_{5}} (
\tilde{v}^{2,2}_{4} x_{4}  x_{4} +
\tilde{v}^{2,2}_{3} x_{3}  x_{5} +
\tilde{v}^{2,2}_{2} x_{2}  x_{6} +
\tilde{v}^{2,2}_{1} x_{1}  x_{1}) 
\end{eqnarray*}
\begin{eqnarray*}
\hat{\mathit{f}}^{3}_{2} &=&  x_{6}  x_{3}+ \frac{1}{\tilde{v}^{3,2}_{5}} (
\tilde{v}^{3,2}_{4} x_{4}  x_{5} +
\tilde{v}^{3,2}_{3} x_{3}  x_{6} +
\tilde{v}^{3,2}_{2} x_{2}  x_{1} +
\tilde{v}^{3,2}_{1} x_{1}  x_{2}) 
\end{eqnarray*}
\begin{eqnarray*}
\hat{\mathit{f}}^{4}_{2} &=&  x_{6}  x_{4}+ \frac{1}{\tilde{v}^{4,2}_{5}} (
\tilde{v}^{4,2}_{4} x_{4}  x_{2} +
\tilde{v}^{4,2}_{3} x_{3}  x_{3} +
\tilde{v}^{4,2}_{2} x_{2}  x_{4} +
\tilde{v}^{4,2}_{1} x_{1}  x_{5}) 
\end{eqnarray*}
\begin{eqnarray*}
\hat{\mathit{f}}^{5}_{2} &=&  x_{6}  x_{5}+ \frac{1}{\tilde{v}^{5,2}_{5}} (
\tilde{v}^{5,2}_{4} x_{4}  x_{1} +
\tilde{v}^{5,2}_{3} x_{3}  x_{2} +
\tilde{v}^{5,2}_{2} x_{2}  x_{3} +
\tilde{v}^{5,2}_{1} x_{1}  x_{4}) 
\end{eqnarray*}
\begin{eqnarray*}
\hat{\mathit{f}}^{6}_{2} &=&  x_{6}  x_{6}+ \frac{1}{\tilde{v}^{6,2}_{5}} (
\tilde{v}^{6,2}_{4} x_{4}  x_{2} +
\tilde{v}^{6,2}_{3} x_{3}  x_{3} +
\tilde{v}^{6,2}_{2} x_{2}  x_{4} +
\tilde{v}^{6,2}_{1} x_{1}  x_{5}) 
\end{eqnarray*}

From the relations one immediately sees that 
\begin{eqnarray*}
\mathbf{S}_{2}' &=& \{  x_{i_1} x_{i_2}  \in  \mathbf{S}_{2} \; | \; i_1 \leq 4 \} 
\end{eqnarray*} 
spans $\mathbb{K}\langle x_1, ..., x_c \rangle _{2}$ modulo $\mathcal{R}$. Equivalently $\mathbf{S}_{2}'$ is a linear basis for $B_g(\tau,\theta)_{2}$. Also, from the discussion above we see that 
\begin{eqnarray*}
\mathbf{S}_{3}' &=& \{  x_{i_1} x_{i_2}  x_{i_3}  \in  \mathbf{S}_{3} \; | \; i_1,i_2 \leq 4 \} 
\end{eqnarray*} 
spans $\mathbb{K}\langle x_1, ..., x_c \rangle_3$ modulo $\mathcal{R}$. $| \mathbf{S}_{3}' |= (a+d)^2c=96$ while $B_g(\tau,\theta)_{3}$ has dimension $((a+d)^2-1)c =90$ thus there are $6$ redundant elements in $\mathbf{S}_{3}'$.

For each one of the relations we look at the obstructions coming from lower relations:
\begin{enumerate}
\item For $\hat{\mathit{f}}^{5}_{1}$:
\begin{eqnarray*}
T(\hat{\mathit{f}}^{5}_{1}) x_{1} = x_{1} T(\hat{\mathit{f}}^{1}_{1}) = x_{5}x_{5}x_{1} &;&
T(\hat{\mathit{f}}^{5}_{1}) x_{2} = x_{5} T(\hat{\mathit{f}}^{2}_{1}) = x_{5}x_{5}x_{2} \\
T(\hat{\mathit{f}}^{5}_{1}) x_{3} = x_{5} T(\hat{\mathit{f}}^{3}_{1}) = x_{5}x_{5}x_{3} &;&
T(\hat{\mathit{f}}^{5}_{1}) x_{4} = x_{5} T(\hat{\mathit{f}}^{4}_{1}) = x_{5}x_{5}x_{3} \\
T(\hat{\mathit{f}}^{5}_{1}) x_{5} = x_{5} T(\hat{\mathit{f}}^{5}_{1}) = x_{5}x_{5}x_{5} 
\end{eqnarray*}
\item For $\hat{\mathit{f}}^{6}_{1}$:
\begin{eqnarray*}
T(\hat{\mathit{f}}^{5}_{1}) x_{6} = x_{5} T(\hat{\mathit{f}}^{6}_{1}) &=& x_{5}x_{5}x_{6} \\
\end{eqnarray*}
\item For $\hat{\mathit{f}}^{1}_{2}$:
\begin{eqnarray*}
T(\hat{\mathit{f}}^{5}_{1}) x_{1} &= x_{5} T(\hat{\mathit{f}}^{1}_{1}) &= x_{5}x_{6}x_{1} \\
\end{eqnarray*}
\item For $\hat{\mathit{f}}^{2}_{2}$:
\begin{eqnarray*}
T(\hat{\mathit{f}}^{5}_{1}) x_{2} &= x_{5} T(\hat{\mathit{f}}^{2}_{2}) &= x_{5}x_{6}x_{2} \\
\end{eqnarray*}
\item For $\hat{\mathit{f}}^{3}_{2}$:
\begin{eqnarray*}
T(\hat{\mathit{f}}^{5}_{1}) x_{3} &= x_{5} T(\hat{\mathit{f}}^{3}_{2}) &= x_{5}x_{6}x_{3} \\
\end{eqnarray*}
\item For $\hat{\mathit{f}}^{4}_{2}$:
\begin{eqnarray*}
T(\hat{\mathit{f}}^{5}_{1}) x_{4} &= x_{5} T(\hat{\mathit{f}}^{4}_{2}) &= x_{5}x_{6}x_{4} \\
\end{eqnarray*}
\item For $\hat{\mathit{f}}^{5}_{2}$:
\begin{eqnarray*}
T(\hat{\mathit{f}}^{5}_{2}) x_{1} = x_{6} T(\hat{\mathit{f}}^{1}_{1}) = x_{6}x_{5}x_{1} &;&
T(\hat{\mathit{f}}^{5}_{2}) x_{2} = x_{6} T(\hat{\mathit{f}}^{2}_{1}) = x_{6}x_{5}x_{2} \\
T(\hat{\mathit{f}}^{5}_{2}) x_{3} = x_{6} T(\hat{\mathit{f}}^{3}_{1}) = x_{6}x_{5}x_{3} &;&
T(\hat{\mathit{f}}^{5}_{2}) x_{4} = x_{6} T(\hat{\mathit{f}}^{4}_{1}) = x_{6}x_{5}x_{4} \\
T(\hat{\mathit{f}}^{5}_{2}) x_{5} = x_{6} T(\hat{\mathit{f}}^{5}_{1}) = x_{6}x_{5}x_{5} &;&
T(\hat{\mathit{f}}^{5}_{2}) x_{6} = x_{6} T(\hat{\mathit{f}}^{6}_{1}) = x_{6}x_{5}x_{6} \\
T(\hat{\mathit{f}}^{6}_{1}) x_{5} = x_{5} T(\hat{\mathit{f}}^{5}_{2}) = x_{5}x_{6}x_{5}.
\end{eqnarray*}
\item For $\hat{\mathit{f}}^{6}_{2}$:
\begin{eqnarray*}
T(\hat{\mathit{f}}^{6}_{2}) x_{1} = x_{6} T(\hat{\mathit{f}}^{1}_{2}) = x_{6}x_{6}x_{1} &;&
T(\hat{\mathit{f}}^{6}_{2}) x_{2} = x_{6} T(\hat{\mathit{f}}^{2}_{2}) = x_{6}x_{6}x_{2} \\
T(\hat{\mathit{f}}^{6}_{2}) x_{3} = x_{6} T(\hat{\mathit{f}}^{3}_{2}) = x_{6}x_{6}x_{3} &;&
T(\hat{\mathit{f}}^{6}_{2}) x_{4} = x_{6} T(\hat{\mathit{f}}^{4}_{2}) = x_{6}x_{6}x_{4} \\
T(\hat{\mathit{f}}^{6}_{2}) x_{5} = x_{6} T(\hat{\mathit{f}}^{5}_{2}) = x_{6}x_{6}x_{5} &;&
T(\hat{\mathit{f}}^{6}_{2}) x_{6} = x_{6} T(\hat{\mathit{f}}^{6}_{2}) = x_{6}x_{6}x_{6} \\
T(\hat{\mathit{f}}^{6}_{1}) x_{5} = x_{5} T(\hat{\mathit{f}}^{6}_{2}) = x_{5}x_{6}x_{6}.
\end{eqnarray*}
\end{enumerate}

For each obstruction of $\hat{\mathit{f}}^{\mu}_{k}$; $T(\hat{\mathit{f}}^{\mu}_{k}) x_{i_3} = x_{i_1} T(\hat{\mathit{f}}^{\mu'}_{k'}) = x_{i_1}x_{i_2}x_{i_3} \in \mathbf{S}_3$ we must look at the normal form of:
$$
\hat{\mathit{f}}^{\mu}_{k} x_{i_3} - x_{i_1} \hat{\mathit{f}}^{\mu'}_{k'} 
$$
with respect to the ideal generated by the relations which are lower than $\hat{\mathit{f}}^{\mu}_{k}$. The following elements give nontrivial normal forms:
\begin{eqnarray*}
\hat{\mathit{f}}^{5}_{2} x_{1} - x_{6} \hat{\mathit{f}}^{1}_{1} &;&
\hat{\mathit{f}}^{5}_{2} x_{2} -x_{6} \hat{\mathit{f}}^{2}_{1}  \\
\hat{\mathit{f}}^{5}_{2} x_{3} -x_{6} \hat{\mathit{f}}^{3}_{1}  &;&
\hat{\mathit{f}}^{5}_{2} x_{4} -x_{6} \hat{\mathit{f}}^{4}_{1}  \\
\hat{\mathit{f}}^{5}_{2} x_{5} - x_{6} \hat{\mathit{f}}^{5}_{1} &;&
\hat{\mathit{f}}^{5}_{2} x_{6} - x_{6} \hat{\mathit{f}}^{6}_{1}.\\
\end{eqnarray*}

The leading terms of the corresponding normal forms must be removed from $\mathbf{S}_{3}'$. These terms are:
\begin{eqnarray*}
x_{4}x_{1}x_{1} ,&   x_{4}x_{1}x_{2} ,&   x_{4}x_{1}x_{3} , \\  
x_{4}x_{1}x_{4} ,&   x_{4}x_{1}x_{5} ,&   x_{4}x_{1}x_{6}.  
\end{eqnarray*} 
Thus a linear basis for $B_g(\tau,\theta)_{3}$ is given by:
\begin{eqnarray}
\{  x_{i_1} x_{i_2}  x_{i_3}  \in  \mathbf{S}_{3} \; | \; i_1,i_2 \leq 4, \; x_{i_1} x_{i_2}  \neq x_{1} x_{4}  \} .
\end{eqnarray} 
\end{example}


\section{Modularity of $B_g(\tau,\theta)$ }
\label{sec7}
Let $g$ satisfy \ref{condiciones1} and \ref{condiciones2} and let $\theta$ be the corresponding real fixed point as given at the end of Section~\ref{sec1}. In this section we consider the algebras $B_g(\tau,\theta)$ as a family of algebras parametrized by $\tau$. We use the fact that the relations defining $B_g(\tau,\theta)$ are given in terms of theta constants which are modular forms of weight $\frac{1}{2}$ and certain level. From this we can get a presentation in which the defining relations have coefficients which are modular functions of the complex structure $\tau$.

Lets start by considering the presentation of $B_g(\tau,\theta)$ given in Theorem~\ref{precentacion1}. Each one of the coefficients $v^{\mu,k}_i$ in the defining relations $\mathit{f}^{\mu}_{k}$  is given by one of the functions $F^{g, \mu}_{i_1,i_2,...,i_{a+d}}(\tau)$. The function $F^{g, \mu}_{i_1,i_2,...,i_{a+d}}(\tau)$ was defined as a determinant of a matrix consisting of theta constants in $\tau'=-l \tau$ therefore it belongs to the ring generated by this functions. Being more precise, let $l$ be a positive integer and let $\mathfrak{C}_l$ be the ring  generated over $\mathbb{C}$ by the functions $\vartheta_{r,s}$ with $(r,s) \in (\frac{1}{l}\mathbb{Z})^2$. Since $\vartheta_{r,s}$ depends on the characteristics $(r,s)$ only modulo integers we see that $\mathfrak{C}_l$ is of finite type over $\mathbb{C}$; it is generated by $\{ \vartheta_{r_i,s_i} \}$ where $(r_i,s_i)$ runs though a set of representatives 
of $(\frac{1}{r}\mathbb{Z}/ \mathbb{Z})^2$ in  $(\frac{1}{r}\mathbb{Z})^2$. If we take $l=c(a+d)$ and view the coefficients as functions of $\tau'=-l \tau$ then our algebras will be naturally defined over $\mathfrak{C}_l$. The ring $\mathfrak{C}_l$  becomes a graded ring by assigning to each $\vartheta_{r,s}\in \mathfrak{C}_l$ weight $\frac{1}{2}$.

The rings $\mathfrak{C}_l$ of theta constants were studied by Igusa in \cite{Igusa1} and \cite{Igusa2}. The main results relate these rings to rings of modular forms for even values of $l$. 

If we denote by $\Gamma_N$ the principal congruence subgroup of level $N$ in $Sl_2(\mathbb{Z})$ then we see from \ref{Gamma2n} that 
\begin{eqnarray}
\Gamma_{2n} \subset \Gamma_{n,2n} \subset \Gamma_{n}
\end{eqnarray}
In particular $\Gamma_{n,2n}$ is a congruence subgroup of $Sl_2(\mathbb{Z})$.

Let $\mathfrak{B}_l:=\mathfrak{C}_l^{(2)}$ be the subring of elements with homogeneous components of even degree i.e. the ring generated over $\mathbb{C}$ by the double products $\vartheta_{r,s} \vartheta_{r',s'}$ with $(r,s), (r',s')$ in $(\frac{1}{l}\mathbb{Z})^2$. The transformation law for theta constants \ref{modulareq} shows that $Sl_2(\mathbb{Z})$ acts on  $\mathfrak{B}_l$ thus one gets a homomorphism from $Sl_2(\mathbb{Z})$ to the group of degree preserving automorphisms of the algebra $\mathfrak{B}_l$. If $l$ is even the kernel of this morphism is precisely $\Gamma_{l,2l^{2}}$. Thus, the normal subgroup of $Sl_2(\mathbb{Z})$ consisting of elements which keep $\mathfrak{B}_l$ element wise invariant is $\Gamma_{l,2l^{2}}$. It follows then that $\tau' \mapsto \vartheta_{r,s}(\tau') \vartheta_{r',s'}(\tau')$ is a modular form of weight $1$ and level $\Gamma_{l^2,2 l^2}$.  We estate this classical result together with a theorem due to Igusa (\cite{Igusa1, Igusa2}). 
\begin{theorem} 
\label{Igusateo}
Let $l$ be a positive even integer. Given a congruence subgroup $\Gamma$ of $Sl_2(\mathbb{Z})$ denote by $\mathfrak{G}(\Gamma)$ the ring of all holomorphic modular forms of level $\Gamma$.  Then:
\begin{enumerate}
\item $\mathfrak{B}_l$ is a subring of $\mathfrak{G}(\Gamma_{l,2l^{2}})$.
\item The integral closure of $\mathfrak{B}_l$ in its field of fractions is $\mathfrak{G}(\Gamma_{l,2l^{2}})$.
\end{enumerate}
\end{theorem}

In order to study the behavior of the coefficients of the defining relations of $B_g(\tau,\theta)$ it is useful at this point to make some remarks about the structure of $\Gamma_{l,2l^{2}}$. They follow from the general theory of discrete subgroups of $Gl_2^{+}(\mathbb{R})$ (c.f. \cite{Shimura}). As above we will consider the action of subgroups of $Gl_2(\mathbb{C})$ on $\mathbb{C}\cup \{\infty \}$ by fractional linear transformations. In particular $Gl_2^{+}(\mathbb{R})$ is identified with the group of holomorphic automorphisms of the upper half plane $\mathbb{H} = \{ \tau \in \mathbb{C} | Im(\tau )>0\}$.

Two subgroups of a group $G$ are said to be commensurable if their intersection has finite index in both of them. If two discrete subgroups $\Gamma$ and $\Gamma '$ of $Sl_2(\mathbb{R})$ are commensurable then their sets of cusps in $\mathbb{C}\cup \{\infty \}$ are the same i.e. the set of points which are fixed points of some parabolic transformation in $\Gamma$ coincides with the corresponding set for $\Gamma '$. In Particular, since $\Gamma_{n,2n}$ is a congruence subgroup of $Sl_2(\mathbb{Z})$ it is commensurable with  it and so the set of cusp of $\Gamma_{n,2n}$ is $\mathbb{Q}\cup \{\infty \}$. Denote by $\mathbb{H}^{*}=\mathbb{H}\cup\mathbb{Q}\cup \{\infty \}$ the upper half plane with the cusps added. It follows from commensurability with $Sl_2(\mathbb{Z})$ that the quotient
\begin{eqnarray}
X_{\Gamma_{n,2n}} = \mathbb{H}^{*} / \Gamma_{n,2n}
\end{eqnarray}
is a compact Riemann surface. In the terminology of Shimura, $\Gamma_{n,2n}$ is a Fuchsian 
group of the first kind (c.f. \cite{Shimura}). 
\begin{notation}
On what follows we will change the sign on our complex structure and 
take $\tau$ with $Im(\tau)>0$. $\tau '$ is then given by $\tau ' = l \tau$.
\end{notation}

In order to study the modularity of the defining relations of $B_g(\tau,\theta)$ as functions of $\tau$ we have to take care of the scaling by $l=c(a+d)$. For this we introduce some notation.
\begin{definition}
\label{gamma}
let $m$ be a even positive integer, define 

\begin{eqnarray}
\Gamma_{n,2n}^{[ m]} =  \Gamma_{1,2}\cap \{ \left(   \begin{array}{cc}
  m & 0 \\
  0 & 1 \\
\end{array} \right)^{-1} \Gamma_{n,2n}   \left( \begin{array}{cc}
  m & 0 \\
  0 & 1 \\
\end{array} \right)  \}
\end{eqnarray}
\end{definition}

This subgroups are the levels for our relations:

\begin{theorem}
\label{precentacion3}
Let  $g$ satisfy \ref{condiciones1} and \ref{condiciones2} and assume $l=c(a+d)$ is even. 
Denote by $w=\lfloor \frac{a+d+1}{2} \rfloor$ the integer part of $\frac{a+d+1}{2}$. Let $\mu \in \{ 1, 2, ...,c \}$,  $k \in \{ 1,..., c-a-d \}$ and let $\tau \mapsto v^{\mu,k}(\tau)$ be given as as in Theorem~\ref{precentacion1}. Define $\hat{v}^{\mu,k} = \hat{v}^{\mu,k}(\tau)\in \mathbb{C}^{c}$ by 
\begin{eqnarray*}
\hat{v}^{\mu,k}(\tau) =  v^{\mu,k}(\tau)  & \text{if $a+d$ is even} \\
\hat{v}^{\mu,k}(\tau) =  \vartheta(l \tau) v^{\mu,k}(\tau)  &  \text{if $a+d$ is odd}.
\end{eqnarray*}
Then:
\begin{enumerate}
\item The algebra $B_g(\tau,\theta)$ is generated by elements $x_1,...,x_c$ of degree $1$ subject to relations $\hat{\mathit{f}}^{\mu}_{k}=0$ where:
\begin{eqnarray*}
\hat{\mathit{f}}^{\mu}_{k}=\hat{v}^{\mu,k}_1 x_{\alpha(\mu,1)} x_1 +...+\hat{v}^{\mu,k}_c x_{\alpha(\mu,c)} x_c
\end{eqnarray*}

\item Each one of the functions $\tau \mapsto \hat{v}^{\mu,k}_j(\tau)$ is a modular form of weight $w$ and level $\Gamma_{l^{2},2l^{2}}^{[ l]}$.
\end{enumerate}
\end{theorem}
\begin{proof}
The first part of the theorem follows from Theorem~\ref{precentacion1}.

For the second part note that each  $\hat{v}^{\mu,k}_j(\tau)$ has the form:  
\begin{eqnarray}
\label{lah}
h^{\mu,k}_j(\tau')=\vartheta(\tau')^{\epsilon}  \sum_{\sigma \in S_{a+d}} \text{sgn} 
(\sigma) \prod_{k=1}^{a+d}
\vartheta_{q(\mu)-\frac{d i_{k}}{c}-\frac{\sigma (k)}{(a+d)},0} (\tau ') 
\end{eqnarray} 
for some $i_1,...,i_{a+d}\in \{ 1,...,a+d \}$ where $\epsilon = 0,1$.

The factor on the right just makes the number of theta constants in the products even so, as function of $\tau'=l\tau$, $h^{\mu,k}_j$ is a homogeneous element of $\mathfrak{B}_l$ of weight $w$. In particular, $h^{\mu,k}_j$ is a  modular form of weight $w$ and level $\Gamma_{l,2l^{2}}$. Since $\hat{v}^{\mu,k}_j(\tau) = h^{\mu,k}_j(l\tau)$ we see that whenever
$$
\left(   \begin{array}{cc}
  x & ly \\
  l^{-1}z & w \\
\end{array} \right) \in \Gamma_{l,2l^{2}}
$$
we have 
$$
\hat{v}^{\mu,k}_j(\frac{x\tau+y}{z\tau +w}) = h^{\mu,k}_j(l \frac{x\tau+y}{z\tau +w}) = (z\tau +w)^k h^{\mu,k}_j(l \tau)
= (z\tau +w)^k \hat{v}^{\mu,k}_j(\tau).
$$
Therefore $\tau \mapsto \hat{v}^{\mu,k}_j(\tau)$ is a modular form of weight $k$ and level 
$$
\Gamma_{1,2}\cap \{ \left(   \begin{array}{cc}
  l & 0 \\
  0 & 1 \\
\end{array} \right)^{-1} \Gamma_{l,2l^{2}}   \left( \begin{array}{cc}
  l & 0 \\
  0 & 1 \\
\end{array} \right)  \}
$$
\end{proof}

The above theorem will allow us to average the values of the coefficients of the defining relations of $B_g(\tau,\theta)$. For this purpose we are
interested in determining whether the modular forms in the defining relations of $B_g(\tau,\theta)$ are modular forms of cusp type; that is, their Fourier expansions at each one of the cusp should have constant term equal to $0$. In what follows we will show that that this is indeed the case.
\begin{proposition}
For $\mu,j \in \{ 1, 2, ...,c \}$ and  $k \in \{ 1,..., c-a-d \}$ let $\hat{v}^{\mu,k}_j(\tau)$ be defined as in Theorem~\ref{precentacion3}. Then $\hat{v}^{\mu,k}_j(\tau)$ is modular form of cusp type for $\Gamma_{l^{2},2l^{2}}^{[ l]}$. 
\end{proposition}
\begin{proof}

By conjugating with  
$\left( \begin{array}{cc}
  l & 0 \\
  0 & 1 \\
\end{array} \right)$ it is enough to show that $h^{\mu,k}_j(\tau')$ in \ref{lah} is a cusp form for $\Gamma_{l,2l^{2}}$. 

We must look at the Fourier series expansions around the cusps of each $h^{\mu,k}_j(\tau')$. Since $h^{\mu,k}_j(\tau')$ is given as a product of theta constants we are interested in the behavior of $\vartheta_{r}$ around the cusps. Moreover, taking into account that the set of cusp forms is an ideal of the graded ring $\mathfrak{G}(\Gamma_{l,2l^{2}})$ of holomorphic modular forms we see that once we show that in each term  of $h^{\mu,k}_j(\tau')$ some factor is a cusp form the result will follow.

Let $l$ be a positive even integer and let $r\in \frac{\mathbb{Z}}{l}$. 
First we look at the behavior of $\vartheta_{r}$ at $\infty$. We can assume $r= \frac{k}{l}$ with $k\in\{ 0,\dots,  l-1 \}$. For $\tau' \in \mathbb{H}$ let $\tilde{q}_{\tau'}=\exp (\frac{ \pi \imath \tau'}{l^{2}})$. Then the series defining $\vartheta_{r}$ is given by: 
\begin{eqnarray*}
\vartheta_{r}(\tau') &=& 
\sum_{n \in \mathbb{Z}} \tilde{q}_{\tau'}^{\phantom{\tau} (nl+rl)^{2}} \\
&=&  \sum_{m \geq 0} a_{m} \tilde{q}_{\tau'}^{\phantom{\tau} m} \\
\end{eqnarray*}
where in the last sum the coefficient $a_{m}$ takes the values $0$, $1$ or $2$ depending on 
whether $m$ is the square of one integer of the form $(nl+rl)^{2}$ for some $n \in \mathbb{Z}$, for two integers of this form or for none. Since we choose $r= \frac{k}{l}$ with $k\in\{ 0,...l-1 \}$ for the constant term we have that $a_{0} \neq 0$ only if $nl+rl=0$ for some $n\in \mathbb{Z}$. This can only happen if $n=r=0$. Therefore the constant term in the $\tilde{q}_{\tau'}$ series expansion of any product of theta constants $\prod_{i=1}^{s} \vartheta_{r_i}$ will vanish provided that at least one of the factors has a nonzero characteristic $r_j \neq 0$. Finally, since any cusp can be carried to $\infty$ by an element in $Sl_2$ the result follows.
\end{proof}

\noindent\begin{example}
\label{superejemplo}
Let  
\begin{eqnarray}
\label{lamatriz}
g &=& \left(
\begin{array}{cc}
  5 & -1 \\
  6 & -1 \\
\end{array}
\right)
\end{eqnarray} 
The eigenvalues of $g$ are 
$$
\lambda^{+} = 2 - {\sqrt{3}}, \quad  \lambda^{-} =  2 + {\sqrt{3}}
$$ 
and the fixed points of $g$ are 
$$
\theta = \frac{3 - {\sqrt{3}}}{6}, \quad  \theta' = \frac{3 + {\sqrt{3}}}{6}.
$$
For any complex structure $\tau$ on $\mathcal{A}_{\theta}$ the corresponding connection $ \bar{\nabla}_0$ on  $E_{-1,6}(\theta)$ will have a six dimensional space of holomorphic sections $\mathcal{H}_g \simeq \mathbb{C}^{6}$. Also,  $l=c(a+d)=24$ so $\mathcal{H}_{g^2}\simeq \mathbb{C}^{24}$. After choosing a basis the multiplication map $m:\mathcal{H}_{g}\otimes \mathcal{H}_{g}\rightarrow \mathcal{H}_{g^2}$ is represented by a $36\times 24$ matrix whose coefficients belong to $\mathfrak{C}_{24}$ when viewed as functions of $\tau'=24\tau$:
$$
C^{\gamma}_{\alpha , \beta} \quad =\quad 
\begin{cases} 
 \quad \vartheta_{\frac{4 \beta - 5 \gamma }{24}} ( \tau' ) 
 \quad if\quad  \alpha  \equiv d(\gamma - \beta) \mod 6 \\
 \quad 0 \quad \text{otherwise}
\end{cases}
$$
We write it as $M\simeq M(1)\oplus  M(2)\oplus  M(3)\oplus  M(4)\oplus  M(5) \oplus  M(6)$ where the elements of $M(\mu )\in  \mathcal{M}_{4,6}(\mathbb{C})$ are given by 
$$
M(\mu )_{i,j}= \vartheta_{\frac{6 - 5 \mu}{24}+\frac{j}{6}-\frac{i}{4}} (\tau' )
$$
The corresponding matrix of characteristics is:
$$
\Lambda = \frac{1}{24}\left(
\begin{matrix}
2 & 22 & 18 & 14 & 10 & 6 \cr 8 & 4 & 0 & 20 & 16 & 12 \cr 14 & 10 & 6 & 2 & 22 & 18 \cr 20 & 16 & 12 & 8 &
   4 & 0 \cr 
\end{matrix}
\right)
$$
Each $\mu\in\{1, ...,6 \}$ gives us a set of $2=c-a-d$ relations corresponding to a basis for the kernel of the matrix $M(\mu)$. Thus $B_{g}(\theta,\tau)$ is in this case is a quadratic algebra with $6$ generators of degree $1$ and $12$ quadratic relations. The corresponding functions of $\tau$ determining the coefficients of the relations are the minors
of each $M(\mu)$. Each ordered 4-tuple $i_1,i_2,i_3, i_4  \in \{ 1,2,3,4,5,6 \}$ determines a weight 2 modular form of $\tau'$ which belongs to $\mathfrak{B}_4$
\begin{eqnarray*}
F^{g, \mu}_{i_1,i_2,i_{3},i_{4}}
& = &  \sum_{\sigma \in S_{4}} \text{sgn} (\sigma) \prod_{k=1}^{4}
\vartheta_{\frac{6 - 5 \mu}{24}-\frac{d i_{k}}{6}-\frac{\sigma (k)}{4}}.  \\
\end{eqnarray*}

Considered as functions of $\tau$ each $F^{g, \mu}_{i_1,i_2,i_{3},i_{4}}$ is a modular form of weight $2$ and level $\Gamma_{24^{2},2(24^{2})}^{[ 24]}$. In particular, each $F^{g, \mu}_{i_1,i_2,i_{3},i_{4}}(\tau)$ defines a differential of the
first kind in the modular curve $X_{\Gamma_{24^{2},2(24^{2})}^{[ 24]}}$.
\end{example}


\section{Modular symbols and averaged algebras}

In this section we use the results about modularity obtained in Section~\ref{sec7} in order to define algebras which do not depend on a particular choice of the complex structure $\tau$. Fix $g$, $l$ and $\theta$ as in the previous sections. Let $\Gamma = \Gamma_{l^{2},2l^{2}}^{[ l]}$, by Theorem~\ref{precentacion3} one can take a presentation of $B_g(\tau,\theta)$ in which each one of the coefficients $v(\tau)$ in the defining relations is a modular form for $\Gamma$. If $v(\tau)$ is a modular form of even weight $w=2r$ we can consider it as a $r$-fold differential on $\mathbb{H}$ invariant under $\Gamma$. That is $v$ is a function in $(dz)^{-r} ((\Omega^{1}_{\mathbb{H}})^{\otimes r})^{\Gamma}$. A $\Gamma$ invariant $k$-fold differential can be pushed down to a differential form in the $(2r-2)$ fibered power of the universal curve over $X_{\Gamma}$. The corresponding integrals along homology classes can be realized as values of line integrals along geodesics in $\mathbb{H}$. This formalism was developed by Manin in \cite{manin4} for modular forms of weight $2$ and extended to arbitrary weights by Shokurov in \cite{Shokurov1}.  

\vspace{0.3 cm} \noindent Let $\Gamma$ be a congruence subgroup of $SL_2 (\mathbb{Z})$.  Let $k$ be a positive integer and consider the action of the crossed product  $\Gamma \ltimes (\mathbb{Z}^{k}\times \mathbb{Z}^{k})$ on $\mathbb{H} \times \mathbb{C}^{k}$ 
given by:
\begin{eqnarray}
(\gamma; n,m): (\tau,z) \mapsto (\gamma \tau, \frac{z + \tau n + m}{C \tau +D} )
\end{eqnarray}
where $\gamma= \left(
\begin{array}{cc}
  A & B \\
  C & D \\
\end{array}
\right)\in \Gamma
; n,m \in \mathbb{Z}^{k}; \tau \in \mathbb{H}$ and $z \in \mathbb{C}^{k}$.

The quotient $\Gamma \ltimes (\mathbb{Z}^{k}\times \mathbb{Z}^{k}) \setminus \mathbb{H} \times \mathbb{C}^{k}$ admits a canonical smooth compactification $\mathcal{Z}_{\Gamma}^{k}$ called the Kuga-Sato variety. Let $\delta, \rho \in \mathbb{P}^{1}(\mathbb{Q})$ be two cusps and let $n,m \in \mathbb{Z}^{k}$. This data defines a homology class $\{ \delta, \rho, n,m\}_{\Gamma} \in H_{k+1}(\mathcal{Z}_{\Gamma}^{k})$ called a modular symbol (c.f. \cite{Shokurov1}). 

\vspace{0.3 cm} \noindent Now let $v(\tau)$ be a cusp modular form of even weight $w=2r$  and level $\Gamma$. Then $\omega = v(\tau)d\tau \wedge d\zeta \wedge d\zeta_1\wedge ... \wedge d\zeta_{w-2}$ is a $\Gamma \times \mathbb{Z}^{w-2} \times \mathbb{Z}^{w-2}$ invariant holomorphic volume form on $\mathbb{H} \times \mathbb{C}^{w-2}$ so it can be pushed down to a holomorphic volume form  $\hat{\omega}$ in the Kuga-Sato variety $\mathcal{Z}_{\Gamma}^{w-2}$. The pairing of this form with the modular symbol $\{ \delta, \rho, n,m \} _{\Gamma}$ is given by (c.f. \cite{manin3}): 
\begin{eqnarray}
\label{modular}
\int_{\delta}^{\rho}v(\tau) \sum_{i=1}^{2r -2} ( n_i \tau + m_i )d\tau 
= \int_{\{ \delta, \rho, n,m \} _{\Gamma}} \hat{\omega}. 
\end{eqnarray}

\vspace{0.3 cm} \noindent Where the integral on the left is the line integral of the holomorphic
differential $v(\tau) \sum_{i=1}^{2r -2} ( n_i \tau + m_i )d\tau$ along the geodesic in $\mathbb{H}$ joining $\delta$ and $\rho$.
\begin{definition}
Let $g$ satisfy (\ref{condiciones1}) and (\ref{condiciones2}) with $l=c(a+d)$   
and $w=\lfloor \frac{a+d+1}{2} \rfloor$ even.  
Let $\hat{v}^{\mu,k}_j(\tau)$ be given as in Theorem~\ref{precentacion3}.
Let $\{ \delta, \rho, n,m \} _{\Gamma} $ be a modular symbol.
We define 
$$
B_g(\theta)\{ \delta, \rho, n,m \} _{\Gamma},
$$ 
the \emph{averaged homogeneous coordinate ring of $\mathbb{T}_{\theta}$ with respect to 
$\{ \delta, \rho, n,m \} _{\Gamma}$} as the quadratic algebra generated by elements $x_1,...,x_c$ of degree $1$ subject to relations:
\begin{eqnarray}
\xi^{\mu}_{k}=\nu^{\mu,k}_1 x_{\alpha(\mu,1)} x_1 +...+\nu^{\mu,k}_c x_{\alpha(\mu,c)} x_c &=&0
\end{eqnarray}
where $\mu \in \{ 1, 2, ...,c \}$, $k \in \{ 1,..., c-a-d \}$ and 
\begin{eqnarray}
\label{integral}
\nu^{\mu,k}_j=\nu^{\mu,k}_j(\{ \delta, \rho, n,m \} _{\Gamma}):= \int_{\delta}^{\rho}\hat{v}^{\mu,k}_j(\tau) \sum_{i=1}^{2r -2} ( n_i \tau + m_i )d\tau
\end{eqnarray}
\end{definition} 
\noindent\begin{example}
\label{superejemplo2}
Lets look at Example~\ref{superejemplo} in this setting.
So we take $$g= \left(
\begin{array}{cc}
  5 & -1 \\
  6 & -1 \\
\end{array}
\right)$$ 
and 
$$
\theta = \frac{3 - {\sqrt{3}}}{6}, 
$$
We are working with modular forms of weight $2$ for the groups $\Gamma_{24^{2},2(24^{2})}^{[ 24]}$ and $\Gamma_{24^{2},2(24^{2})}$. As remarked in Example~\ref{superejemplo}, having weight 2, 
the  modular forms appearing as coefficients of the defining relations of $B_g(\theta, \tau)$
correspond to differentials of the first kind in the modular curve  $X_{\Gamma_{24^{2},2(24^{2})}^{[ 24]}}$. The corresponding modular symbols 
$$
\{ \delta, \rho \} _{\Gamma_{24^{2},2(24^{2})}^{[ 24]}} \in H^{1}(X_{\Gamma_{24^{2},2(24^{2})}^{[ 24]}},\mathbb{Q})
$$
are classical and given any two cusp $\delta, \rho \in \mathbb{Q}\cup \{\infty \}$ we get a averaged homogeneous coordinate ring $B_g(\theta)\{ \delta, \rho \} _{\Gamma_{24^{2},2(24^{2})}^{[ 24]}}$ whose defining relations have coefficients 
$$
\nu^{\mu,k}_j(\{ \delta, \rho\} _{\Gamma_{24^{2},2(24^{2})}^{[ 24]}}):= \int_{\delta}^{\rho}\hat{v}^{\mu,k}_j(\tau) d\tau .
$$
\end{example}


In the definition of the the averaged homogeneous coordinate ring 
$B_g(\theta)\{ \delta, \rho, n,m \} _{\Gamma}$ we  have a non canonical choice coming from the 
cusps $ \delta $ and $\rho$ in the limits of the integration \ref{integral}. We will use the 
ideas developed by Manin and Marcolli in \cite{ManinMarcolli} to get canonical set of cusps 
associated to $\theta$. The modular symbols corresponding to these cusps will then define  
a averaged homogeneous coordinate ring canonically associated to 
$\mathbb{T}_{\theta}$.

One can not naively replace one of the cusp in $\{ \delta, \rho\} _{\Gamma}$  by a irrational number  $\beta \in  \mathbb{R} \setminus \mathbb{Q}$ since corresponding integral \ref{integral} would then diverge. Still, it makes sense to look at the asymptotic behavior of the integrals along infinite geodesics in $\mathbb{H}$ having a irrational endpoint $\beta$. If this asymptotic limit exist it defines a {\it limiting modular symbol} (c.f. \cite{ManinMarcolli}): 
\begin{eqnarray}
 \{ \{ *, \beta \} \}_{\Gamma}  \in H^{1}(X_{\Gamma},\mathbb{R}).
\end{eqnarray}
If $\beta$ is a real quadratic irrationality then the corresponding limiting modular symbol exist and can be computed as a combination of classical modular symbols.
\begin{example}
\label{particular}
Consider again the situation in Example~\ref{superejemplo2}. Thus $\theta= \frac{3 - {\sqrt{3}}}{6}$ and $g$ is given by (\ref{lamatriz}). Let 
$$
\tilde{g}= g^{4 l^{2}} =g^{4 (24^{2})}
$$
Then we have that $\tilde{g}$ is a hyperbolic element of $\Gamma_{24^{2},2(24^{2})}^{[ 24]}$ and $\theta$ is one of its fixed points. Let  also $\tilde{\lambda}^{-} = (\lambda^{-})^{4 l^{2}}>1$ denote the corresponding eigenvalue of $\tilde{g}$. The limiting modular symbol defined by  $\theta$ is given in this case by (c.f. \cite{ManinMarcolli}):
\begin{eqnarray}
\label{limod}
 \{ \{ *, \theta \} \}_{\Gamma_{24^{2},2(24^{2})}^{[ 24]}} = \frac{ \{ 0, \tilde{g}(0) \} _{ \Gamma_{24^{2},2(24^{2})}^{[ 24]}} }{\log \tilde{\lambda}^{-}}
\in H^{1}(X_{\Gamma_{24^{2},2(24^{2})}^{[ 24]}},\mathbb{R}).
\end{eqnarray}

We can now integrate along this homology class the modular forms defining $B_g(\tau,\theta)$. Taking  
$$
\nu^{\mu,k}_j(\theta) := \int_ { \{ \{ *, \theta \} \}_{\Gamma_{24^{2},2(24^{2})}^{[ 24]}}}\hat{v}^{\mu,k}_j(\tau) d\tau
$$
and imposing the relations $\nu^{\mu,k}_1 (\theta) x_{\alpha(\mu,1)} x_1 +...+\nu^{\mu,k}_c (\theta) x_{\alpha(\mu,c)} x_c$ on $ \mathbb{C}\langle x_1, ..., x_c \rangle$ we get a set of relations for a quadratic algebra $B_g(\theta)$ canonically associated with $\theta$ and $g$.
\end{example}

The limiting modular symbol defined by a real  quadratic irrationality $\theta\in (0,1)$ can be computed using its continued fraction expansion. Let $\{ k_n(\theta) \; | \; n=1,2,... \} \subset \mathbb{N}$ be the eventually periodic sequence corresponding to the continued fraction expansion of $\theta$. The corresponding convergents are: 
\begin{eqnarray}
[k_1(\theta),k_2(\theta),...,k_n(\theta)] &= \frac{1}{k_1(\theta)+\frac{1}{k_2(\theta)
+....\frac{1}{k_n(\theta)}}} &= \frac{p_n(\theta)}{q_n(\theta)} .
\end{eqnarray}
\vspace{0.3 cm} \noindent Let also 
\begin{eqnarray}
\label{masmatriz}
g_n(\theta):= \left(
\begin{array}{cc}
  p_{n-1}(\theta) & p_{n}(\theta) \\
  q_{n-1}(\theta) & q_{n}(\theta) \\
\end{array}
\right)\in GL_2 (\mathbb{R})
\end{eqnarray}
and take $\lambda(\theta)= \lim_{n \to \infty} \frac{\log q_n(\theta)}{n}$. Then the modular symbol defined by $\theta$ can be computed by the following formula (c.f.\cite{Marcolli}):
\begin{eqnarray}
\label{cong}
\{ \{ *, \theta \} \}_{\Gamma} &=&
\sum_{n =1}^{m}  \frac{ \{ g_{n}^{-1}(\theta) 0 \, ,\,  g_{n}^{-1}(\theta) \imath \infty \} _{ \Gamma}}{
m \lambda(\theta) }.
\end{eqnarray}
Where $m$ is the period of the continued fraction expansion $\{ k_n(\theta) \; | \; n=1,2,... \} $.

At present the theory of limiting modular symbols has been developed only for weight $w=2$. It is expected that an analogous theory for higher weight can be developed. For our purposes it is enough to consider \ref{cong} as providing a canonical choice of cusps defining  modular symbols over which to average the coefficients of the defining relations of the homogeneous coordinate ring $B_g(\tau,\theta)$:
\begin{definition}
Let  $g$ satisfy (\ref{condiciones1}) and (\ref{condiciones2}) with $l=c(a+d)$ and $w=\lfloor \frac{a+d+1}{2} \rfloor$ even. Take $\Gamma = \Gamma_{l^{2},2l^{2}}^{[ l]}$ and let $\hat{v}^{\mu,k}_j(\tau)$ be given as in Theorem~\ref{precentacion3}. Let also $g_n(\theta), \lambda(\theta)$ and $m$ be as in (\ref{cong}). We define $B_g(\theta)$ the \emph{averaged homogeneous coordinate ring of $\mathbb{T}_{\theta}$}as the quadratic algebra generated by elements $x_1,...,x_c$ of degree $1$ subject to relations:
\begin{eqnarray}
\hat{\xi}^{\mu}_{k}=\hat{\nu}^{\mu,k}_1 x_{\alpha(\mu,1)} x_1 +...+\hat{\nu}^{\mu,k}_c x_{\alpha(\mu,c)} x_c &=& 0 
\end{eqnarray}
where $\mu \in \{ 1, 2, ...,c \}$, $k \in \{ 1,..., c-a-d \}$ and
\begin{eqnarray}
\label{integral2}
\hat{\nu}^{\mu,k}_j &=& \frac{1}{m \lambda(\theta)}
\sum_{n =1}^{m} \int_{g_{n}^{-1}(\theta) 0}^{g_{n}^{-1}(\theta) \imath \infty}\hat{v}^{\mu,k}_j(\tau)d\tau .
\end{eqnarray}
\end{definition} 


\section{The geometric data}

The role played by modular forms in the above discussion points to deep relations with the quantum thermodynamical system introduced by Connes and Marcolli in \cite{ConnesMarcolli} in relation with the class field theory of the modular field. Several results point to the fact that quantum statistical mechanics provides the right framework under which the tools of noncommutative geometry may be  applied to class field theory.

In the case of real quadratic fields explicit class field theory is conjecturally given in terms of special values of $L$-functions, this is the content of Stark's conjectures \cite{Stark}. In order to apply our results on noncommutative tori in this direction using the tools of quantum statistical mechanics we still need to find $C^{*}$-completions of the algebras $B_g(\tau,\theta)$, $B_g(\theta)\{ \delta, \rho, n,m \} _{\Gamma}$ and $B_g(\theta)$. The constructions of \cite{CoDV} seem to provide the right tools in for this purpose. 

The first step in this direction is the construction of the geometric data corresponding to the algebra $B_g(\theta,\tau)$. The geometric data associates to a finitely generated graded algebra $A=\oplus A_n$ a triple $\mathit{T}=(Y,\sigma,\mathcal{L})$  where $Y$ is a projective variety, $\sigma$ is an automorphism of $Y$ and $\mathcal{L}$ is an ample line bundle over $Y$. Starting from such a triple one can construct the graded algebra:
\begin{eqnarray}
B(\mathit{T})= \bigoplus_{n\geq 0} 
H^0(Y,\mathcal{L} \otimes\mathcal{L}^{\sigma}\otimes...\otimes\mathcal{L}^{\sigma^{n-1}}) 
\end{eqnarray}
where $\mathcal{L}^{\sigma} := \sigma^{*} \mathcal{L}$ 
and the multiplication of two sections 
$s_{1}\in B(\mathit{T})_n,s_{2}\in B(\mathit{T})_m$ is given by $s_1s_2:=s_1\otimes s_2^{\sigma^{n}}$.

The construction is made in such a way that one gets a morphism $A\rightarrow B(\mathit{T})$. This construction was introduced by Artin, Tate and Van den Bergh in \cite{ArtinTate} in order to study regular algebras of dimension 3.

Let $\mathcal{T} =\mathbb{K} \langle x_1, ..., x_c \rangle$ be the free associative algebra on $c$ generators of degree one over $\mathbb{K}$. So $\mathcal{T}_1 = \sum \mathbb{K}x_i \simeq \mathbb{K}^{c}$ and $\mathcal{T}$ is the tensor algebra 
\begin{eqnarray}
\mathcal{T} = \bigoplus_{n \geq 0}(\mathcal{T}_1)^{\otimes n} \simeq \bigoplus_{n \geq 0}(\mathbb{K}^{c})^{\otimes n}.
\end{eqnarray}
To each homogeneous element $f\in \mathcal{T}_n=\mathcal{T}_1^{\otimes n}$ we associate the corresponding n-multilinear form $\check{f}: \mathcal{T}_1^{*}\times...\times \mathcal{T}_1^{*}\rightarrow \mathbb{K}$ acting on the n-th Cartesian product of the dual space $\mathcal{T}_1^{*}$. We call $\check{f}$ the multi linearization of $f$. Since $\check{f}$ is multihomogeneous its zero locus defines a hypersurface in $(\mathbb{P}^{c-1}(\mathbb{K}))^{n}$. 

Let now $A$ be a finitely generated quadratic algebra over $\mathbb{K}$.
Assume $A$ is generated in degree one so that 
$$
A \simeq \mathcal{T}/ \mathcal{R}
$$
where $\mathcal{R}=(f_1,\dots,f_r), \; f_i \in \mathcal{T}_1 \otimes \mathcal{T}_1$ is the homogeneous ideal generated by the defining relations of the algebra $A$. The locus of common zeroes of the multilinearizations of the elements of $\mathcal{R}$ defines a variety $\{ \check{f}_i = 0\} = \Gamma\subset \mathbb{P}^{c-1}\times \mathbb{P}^{c-1}$. Let $Y_1$ and $Y_2$ be the corresponding projections and $\sigma: Y_1 \rightarrow Y_2$  be the correspondence with graph $\Gamma$. Assume we can make an identification $Y=Y_1=Y_2$. In the case $\sigma$ is an isomorphism we consider it as an automorphism of $Y$. Letting then 
$i:Y\hookrightarrow \mathbb{P}^{c-1}$ be the inclusion and taking $\mathcal{L}=i^{*}\mathcal{O}_{\mathbb{P}^{c-1}}(1)$ we get the corresponding geometric data $\mathit{T}=(Y,\sigma,\mathcal{L})$. Taking $B(\mathit{T})$ as above we have that the canonical map $A_1 \rightarrow H^0(Y,\mathcal{L})$ extends to a morphism of graded algebras $A \rightarrow B(\mathit{T})$. We call $Y$ the characteristic variety of $A$.

Consider now the algebra $B_g(\theta,\tau)$. We start with its presentation in terms of generators and relations given in Theorem~\ref{precentacion3}

For $\mu \in \{ 1, 2, ...,c \}$ and  $k \in \{ 1,..., c-a-d \}$  let 
\begin{eqnarray}
\mathit{f}^{\mu}_{k}= v^{\mu,k}_1 x_{\alpha(\mu,1)} x_1 +...+v^{\mu,k}_c x_{\alpha(\mu,c)} x_c
\end{eqnarray}
be the corresponding quadratic relation. Denote by $(x_i)_1(x_j)_2$ the map $\mathbb{C}^{c}\times \mathbb{C}^{c} \rightarrow \mathbb{C}$ given by $(v,w)\mapsto v_i w_j$. The multilinearization of $\mathit{f}^{\mu}_{k}$ is then 
\begin{eqnarray}
\check{\mathit{f}}^{\mu}_{k}= v^{\mu,k}_1 (x_{\alpha(\mu,1)})_1 (x_1)_2 +...
+v^{\mu,k}_c (x_{\alpha(\mu,c)})_1 (x_c)_2
\end{eqnarray}

By definition the graph $\Gamma\subset \mathbb{P}^{c-1}\times \mathbb{P}^{c-1}$ of the correspondence $\sigma$ in the geometric data of $B_g(\theta,\tau)$ is the common zero locus of the $c(c-a-d)$ bihomogeneous forms $\check{\mathit{f}}^{\mu}_{k}$. Let $\Omega \in \mathcal{M}_{c,c(c-a-d)}(\mathcal{T}_1)$ be the matrix defined by 
\begin{eqnarray}
\mathit{f}^{\mu}_{k}=\Omega^{\mu}_{k,i}x_i, &  i=1,...,c; \, \mu =1, 2, ...,c;\, k = 1,..., c-a-d.
\end{eqnarray}
If the images of $\Gamma$ under the two projections 
$\pi_1: \mathbb{P}^{c-1}\times \mathbb{P}^{c-1} \rightarrow \mathbb{P}^{c-1}$ and 
$\pi_2: \mathbb{P}^{c-1}\times \mathbb{P}^{c-1} \rightarrow \mathbb{P}^{c-1}$ are equal then $\sigma$ is an automorphism of $Y=\pi_1(\Gamma)=\pi_2(\Gamma)$. $Y$ is given in this case by the vanishing of the $c\times c$ minor determinants of the matrix: 
\begin{eqnarray*}
\check{\Omega} = \left(
\begin{array}{cccc}
v^{1,1}_1 (x_{\alpha(1,1)})_1  &  v^{1,1}_2 (x_{\alpha(1,2)})_1 &\quad ...\quad & v^{1,1}_c (x_{\alpha(1,c)})_1  \\
v^{1,2}_1 (x_{\alpha(1,1)})_1  &  v^{1,2}_2 (x_{\alpha(1,2)})_1 &\quad ...\quad & v^{1,2}_c (x_{\alpha(1,c)})_1  \\
\vdots & \vdots &\quad ...\quad & \vdots  \\
v^{1,c-a-d}_1 (x_{\alpha(1,1)})_1  &  v^{1,c-a-d}_2 (x_{\alpha(1,2)})_1 &\quad ...\quad & v^{1,c-a-d}_c (x_{\alpha(1,c)})_1  \\
v^{2,1}_1 (x_{\alpha(2,1)})_1  &  v^{2,1}_2 (x_{\alpha(2,2)})_1 &\quad ...\quad & v^{2,1}_c (x_{\alpha(2,c)})_1  \\
\vdots & \vdots &\quad ...\quad & \vdots  \\
v^{2,c-a-d}_1 (x_{\alpha(2,1)})_1  &  v^{2,c-a-d}_2 (x_{\alpha(2,2)})_1 &\quad ...\quad & v^{2,c-a-d}_c (x_{\alpha(2,c)})_1  \\
\vdots & \vdots &\quad ...\quad & \vdots  \\ \vdots & \vdots &\quad ...\quad & \vdots  \\
v^{c,1}_1 (x_{\alpha(c,1)})_1  &  v^{c,1}_2 (x_{\alpha(c,2)})_1 &\quad ...\quad & v^{c,1}_c (x_{\alpha(c,c)})_1  \\
\vdots & \vdots &\quad ...\quad & \vdots  \\
v^{c,c-a-d}_1 (x_{\alpha(c,1)})_1  &  v^{c,c-a-d}_2 (x_{\alpha(c,2)})_1 &\quad ...\quad & v^{c,c-a-d}_c (x_{\alpha(c,c)})_1  
\end{array}
\right).
\end{eqnarray*}

\vspace{0.3 cm} \noindent From this we see that $Y\subset  \mathbb{P}^{c-1}$ is defined by ${c(c-a-d) \choose c}$ homogeneous equations of degree $c$. Each one of these equations has as coefficients values of theta constants in $\mathfrak{B}_l$. The explicit form of these determinantial varieties depends in each particular case on the relations satisfied by the corresponding theta constants. 
\begin{remark}
In \cite{Smith} Smith and Stafford analyzed the geometric data of a class of graded rings known as Sklyanin algebras. This are quadratic algebras whose relations have coefficients which can be realized as values of theta functions with characteristics in $\frac{1}{2}\mathbb{Z}$. 
The relations satisfied by $\vartheta_{0,0},\vartheta_{0,\frac{1}{2}},\vartheta_{\frac{1}{2},0}$ and $\vartheta_{\frac{1}{2},\frac{1}{2}}$ play a central role on their analysis and are responsible for the appearance of an elliptic curve as part of the characteristic variety. At a firsts glance it seems like our situation is analogous to that on \cite{Smith} and either
the modular curves or $Proj \mathfrak{B}_l$ could be playing the role that elliptic curves played for the Sklyanin algebras. 
\end{remark}

We shall return to this point in future work where a detailed analysis of the geometric data associated to the various rings considered in this article will be carried out. We expect then the use of the techniques  developed in \cite{CoDV} will then make it possible to  obtain $C^{*}$-completions suitable for arithmetic applications.
\vspace{0.3 cm} 

\appendix
\section{Theta functions and theta constants with rational characteristics}
\label{appendix}
In this appendix we recall the main facts about theta functions and theta constants with rational characteristics that are used in this article. Our treatment follows closely \cite{Mumford2}, \cite{tata1} and \cite{tata3}.

Let $(r,s)\in \mathbb{Q}^2$. The series 
\begin{eqnarray}
\label{thetaseries}
\vartheta_{r,s}(z,\tau) = \sum_{n \in \mathbb{Z}} \exp [\pi \imath (n+r)^{2} \tau +2 \pi \imath (n+r)(z+s) ]
\end{eqnarray}
defines a holomorphic function of $(z,\tau)\in \mathbb{C} \times \mathbb{H}$.  We call 
$\vartheta_{r,s}(z,\tau)$ the {\it theta function with rational characteristics $(r,s)$.} Note that changing the characteristics by integer values does not affect the values of $\vartheta_{r,s}(z,\tau)$. Taking $(r,s)=(0,0)$ we get Riemann's theta function $\vartheta(z,\tau):= \vartheta_{0,0}(z,\tau)$. For $(r,s)\in \mathbb{Q}^2$ we have 
\begin{eqnarray}
\label{thetaseries2}
\vartheta_{r,s}(z,\tau)= \exp [\pi \imath r^{2}\tau +2 \pi \imath r(z+s) ] \vartheta(z+r\tau +s,\tau).
\end{eqnarray}

One of the most important and useful result about $\vartheta$ is the functional equation it satisfies:
\begin{theorem}
\label{funceq}
Let 
$
\gamma= \left(
\begin{array}{cc}
  A & B \\
  C & D \\
\end{array}
\right)\in \Gamma_{1,2}$. Then
\begin{eqnarray*}
 \vartheta ( \frac{z}{C \tau +D},\frac{A \tau +D}{C \tau +D} )= 
\kappa(\gamma) ({C \tau +D})^{\frac{1}{2}} \exp [\frac{\pi \imath C z^{2}}{C \tau +D} ] \vartheta(z, \tau)
\end{eqnarray*}
where the constant $\kappa(\gamma)$ satisfies $\kappa(\gamma)^{8}=1$ and 
$\Gamma_{n,2n}$ is defined for any positive integer $n$ by 
\begin{eqnarray}
\Gamma_{n,2n} &=& \{ \gamma = \left( 
\nonumber
\begin{array}{cc}
  x & y \\
  z & w \\
\end{array}
\right)\in SL_2 (\mathbb{Z})
\, |\,  
 \gamma \equiv 1_{2}\mod n ; xy \equiv  zw \equiv 0 \mod 2n 
 \} \\
\label{Gamma2n}
\end{eqnarray}
\end{theorem}

For a fixed $\tau \in  \mathbb{H}$ consider the lattice $\Lambda_{\tau}= \mathbb{Z} \oplus \tau \mathbb{Z}$. Let $(r,s)\in \mathbb{Q}^2$ and define for 
$\lambda = m' + m \tau \in \Lambda_{\tau}$ 
\begin{eqnarray*}
e_{\lambda}(z):= \exp [- \pi \imath m^{2} \tau +2 \pi \imath (r m' -m (z+s)]
\end{eqnarray*}
Then 
\begin{eqnarray*}
\vartheta_{r,s}(z+ \lambda ,\tau)= e_{\lambda} \vartheta_{r,s}(z,\tau)
\end{eqnarray*}
and the values $e_{\lambda}$ satisfy:
\begin{eqnarray*}
e_{\lambda_1 +\lambda_2}(z)= e_{\lambda_1}(z+ \lambda_2) e_{\lambda_2}(z).
\end{eqnarray*}
Thus the function $\vartheta_{r,s}(z,\tau)$ is an entire $\Lambda_{\tau}$-automorphic function of $z \in \mathbb{C}$. The functions $\{ e_{\lambda} \}_{\lambda \in \Lambda_{\tau}}$ are called automorphy factors. The above conditions just mean that each 
$\vartheta_{r,s}(z,\tau)$ gives a section of a holomorphic line bundle on the elliptic curve $E_{\tau}= \mathbb{C} / \Lambda_{\tau}$. 

We will denote by $\mathcal{L}$ the \emph{basic line bundle over $E_{\tau}$}. This is the bundle obtained as the quotient of $ \mathbb{C}\times  \mathbb{C}$ by the action of $\Lambda_{\tau}$ given by 
\begin{eqnarray*}
\lambda \, \cdot \,  (\zeta, z ) = (\exp [- \pi \imath m^{2} \tau +2 \pi \imath m (z- m')] \zeta ,z- \lambda)
\end{eqnarray*}
where $\lambda = m' + m \tau \in \Lambda_{\tau}, \; (\zeta, z )\in \mathbb{C}\times  \mathbb{C}$, and we view $\mathbb{C}\times  \mathbb{C}$ as the total space of the trivial line bundle over 
$\mathbb{C}$. For an integer $l>1$ the functions $z \mapsto \vartheta_{r,s}(l z,\tau)$ with  $(r,s)\in \frac{1}{l}\mathbb{Z}$ give sections of a line bundle whose $l$ tensor power is isomorphic to 
$\mathcal{L}$. In this case the Lefschetz embedding theorem gives us an embedding to projective space:
\begin{theorem}
\label{Lef}
Let $l>1$ be an integer and let $(r_i,s_i) \in (\frac{1}{l}\mathbb{Z})^2$, $i=0,...,l^2-1 $, be a complete set of representatives of 
$(\frac{1}{l}\mathbb{Z}/ \mathbb{Z})^2$.
Then the holomorphic map $\phi_l:  \mathbb{C} / \Lambda_{\tau} \rightarrow  \mathbb{P}^{l^2 -1}\mathbb{C}$ given by 
\begin{eqnarray}
\phi_l: \quad z \mapsto ( \vartheta_{r_0,s_0}(lz, \tau) ,..., 
\vartheta_{r_{l^2-1},s_{l^2-1}}(l z, \tau) )
\end{eqnarray}
is an embedding. In particular, $\phi_l (\mathbb{C} / \Lambda_{\tau})$ is an algebraic subvariety of $\mathbb{P}^{l^2 -1}\mathbb{C}$. 
\end{theorem}

\begin{remark}
The algebraic equations defining $\phi_l (\mathbb{C} / \Lambda_{\tau})$ can be determined in some cases by relations satisfied by the $\vartheta_{r,s}$. For instance, Riemann's 16 theta relations 
between $\vartheta_{0,0} ,\vartheta_{0,\frac{1}{2}}, \vartheta_{\frac{1}{2},0}$ and 
$\vartheta_{\frac{1}{2},\frac{1}{2}}$ can be used to realize 
$\phi_2 (\mathbb{C} / \Lambda_{\tau})$ as the curve in $\mathbb{P}^{3}\mathbb{C}$ defined by the quadratic equations
\begin{eqnarray*}
\vartheta_{0,0}(0)^{2} X_0^{2}  &=& \vartheta_{0,\frac{1}{2}}(0)^{2} X_1^{2} + \vartheta_{\frac{1}{2},0}(0)^{2} X_2^{2} \cr
\vartheta_{0,0}(0)^{2} X_3^{2}  &=& \vartheta_{\frac{1}{2},0}(0)^{2} X_1^{2} - \vartheta_{0,\frac{1}{2}}(0)^{2} X_2^{2}
\end{eqnarray*}
\end{remark}

The zeroes of $z \mapsto \vartheta_{r,s}(z,\tau)$ are characterized by the following lemma: 
\begin{lemma}
\label{zerostheta}
Fix $\tau \in  \mathbb{H}$, let $l>1$ be an integer and let $(r,s)\in \mathbb{Q}$. Then, the zeroes of $z \mapsto \vartheta_{r,s}(z,\tau)$ are 
the points $(r + p+ \frac{1}{2})\tau + (s + q+ \frac{1}{2})$ with $p,q\in\mathbb{Z}$ .
\end{lemma}

The holomorphic function on $\mathbb{H}$ obtained when we restrict 
$\vartheta_{r,s}(z,\tau)$ to $z=0$ will be refereed as the {\it theta constant with rational characteristics $(r,s)$}. We will use the following notations:
\begin{eqnarray}
\label{notheta}
\vartheta_{r,s}(\tau) &:=& \vartheta_{r,s}(0, \tau).\\
\vartheta_{r}(\tau) &:=& \vartheta_{r,0}(0, \tau).
\end{eqnarray}

The functional equation in Theorem~\ref{funceq} becomes in this case:
\begin{eqnarray}
\label{modulareq}
\vartheta (\frac{A \tau +D}{C \tau +D} ) &=& 
\kappa(\gamma) ({C \tau +D})^{\frac{1}{2}}  \vartheta(\tau).
\end{eqnarray}

A very important variant of these functions is given by the algebraic theta constants:
\begin{eqnarray}
\label{otratheta}
\vartheta_{r,s}^{\alpha}(\tau) \;:=\; e^{-\pi \imath r s  }\vartheta_{r,s}(0, \tau)&& (r,s)
\in \mathbb{Q}^2.
\end{eqnarray}
The rationality behavior of these functions plays an essential role for us. 
Let $\tau \in \mathbb{H}$. Consider as above the elliptic curve $E_{\tau}$ and the basic line bundle $\mathcal{L}$ on it. For all $a,b \in \mathbb{Q}$ we take 
\begin{eqnarray*}
x(a,b)_{m}&=& \text{Image in } E_{\tau} \text{ of } \frac{1}{m}(a \tau +b)
\end{eqnarray*}
Then we have:

\begin{theorem}\emph{(\cite{tata3} Corollary 5.12)}
\label{teotata3}
Suppose that $a,b \in \frac{1}{l}\mathbb{Z}$ and $k \hookrightarrow \mathbb{C}$ is a subfield such that $E_{\tau}$ and $\mathcal{L}$ can be defined over $k$ and $x(a,b)_{2l}$ is rational over 
k. Then
\begin{eqnarray*}
\frac{\vartheta_{a,b}^{\alpha} (\tau) }{\vartheta_{0,0}^{\alpha}(\tau) } &\in& k.
\end{eqnarray*}
\end{theorem}


\bibliographystyle{amsalpha}

\end{document}